\documentclass[global]{svjour}
\usepackage{amssymb}
\usepackage{amsfonts}
\usepackage{amsmath}

\input{tcilatex}

\begin{document}

\title{Existence of equilibrium for generalized games in choice form and
applications}
\author{Monica Patriche}
\institute{University of Bucharest 
%TCIMACRO{\TeXButton{email}{\email{monica.patriche@yahoo.com}}}%
%BeginExpansion
\email{monica.patriche@yahoo.com}%
%EndExpansion
}
\mail{\\
University of Bucharest, Faculty of Mathematics and Computer Science, 14
Academiei Street, 010014 Bucharest, Romania}
\maketitle

\bigskip Abstract. This paper has two central aims: first, to provide simple
conditions under which the generalized games in choice form and,
consequently, the abstract economies, admit equilibrium; second, to study
the solvability of several types of systems of vector quasi-equilibrium
problems as an application. Our work outlines that there still is much to be
gained from using the results concerning the existence of equilibrium of
games as tools of research for other optimization problems.


\section{Introduction and Preliminaries}

The importance of vector equilibrium problems comes, obviously, from the
great number of recent papers dedicated to the study of the existence of
their solutions. We can provide a short list of references presented in our
bibliography: [1-5], [7], [9], [10], [14-18], [21-25], [29]. From the
scientific point of view, we must state the fact that the vector equilibrium
problems unifies several problems, among which we can mention: vector
variational inequalities, vector complementarity problems and vector
optimizations problems.

The purpose of this paper is to provide new conditions under which the
systems of generalized vector quasi-equilibrium problems have solutions.
Following the approach of Lin, Chen and Ansari \cite{lin2007}, we firstly
prove new results based on the existence of equilibrium for generalized
abstract economies, results which extend and improve a number of theorems
existing in literature.

We recall that the notion of abstract economy was introduced by Shafer and
Sonnenschein in \cite{sh}, in order to generalize Nash's model of
noncooperative game in normal form, defined in \cite{80} and \cite{79}.
Shafer and Sonnenschein's model consists of a finite set of agents. Each
agent $i$ has a constraint correspondence $A_{i}$ and, instead of the
utility function $u_{i}$, considered by Nash, he has a preference
correspondence $P_{i}$. Later, several generalizations have been made. Among
them, there are the generalized abstract economies, which have two
constraint correspondences. In a recent paper (\cite{he}), Herv\'{e}s-Beloso
and Patriche provided relevant examples of abstract economies with two
constraint correspondences. The model of Kim and Tan \cite{kim} also
contains an additional "fuzzy" constraint.

In order to establish the existence of equilibrium for a generalized
abstract economy, we apply a new technique, which is completely different
from the classical techniques used in literature, namely the ones proposed
by Shafer and Sonnenschein in \cite{sh} or by Yannelis and Prahbakar in \cite%
{yan}. The method of the last authors consists of the construction, for each
players, of a certain correspondence which has different values, depending
on the set of points for which the intersection of the values of
correspondences $P_{i}$ and $A_{i}$ are empty or not. The product of all
these correspondences are forced to fulfill the conditions of a known fixed
point theorem. Based on this idea, all following results had similar
assumptions, especially which concerned nonempty convex closed (or open)
values of the correspondences $P_{i}$ and $A_{i}$. Our approach consists of
the construction of an auxiliary generalized game in choice form.\medskip

The model of a generalized game in choice form was introduced by Ferrara and
Stefanescu in \cite{fe}. Ferrara and Stefanescu's paper rediscussed the
model of game in choice form, which was firstly defined by these authors in 
\cite{st} and \cite{SS}.

A generalized game in choice form is the family of the individual strategies
sets, the constraint correspondences and a choice profile. A choice profile
can be expressed as a collection of subsets from the set including all
strategies of the game. The difference, from the classical models of Nash (%
\cite{80},\cite{79}) and Debreu (\cite{d}), consists in the fact that
Ferrara and Stefanescu's new model also takes into account the particular
cases when the players' preferences need not be explicitly represented, or
when the choice of a player need not be the best reply to the strategy
combination of the others.\medskip

The rest of this section recalls some definitions and preliminary results
which will be used in the following sections.

Let $X$, $Y$ be topological spaces and $T:X\rightarrow 2^{Y}$ be a
correspondence. The correspondence $\overline{T}$ is defined by $\overline{T}%
(x):=\{y\in Y:(x,y)\in $cl$_{X\times Y}$ Gr $T\}$ (the set cl$_{X\times Y}$
Gr $(T)$ is called the adherence of the graph of $T$)$.$ It is easy to see
that cl $T(x)\subset \overline{T}(x)$ for each $x\in X.$

For each $x\in X,$ the set $T(x)$ is called the \textit{upper section} of $T$
at $x.$ For each $y\in Y,$ the set $T^{-1}(y):=\{x\in X:y\in T(x)\}$ is
called the \textit{lower section} of $T$ at $y$. $T$ is said to have \textit{%
open lower sections} if $T^{-1}(y)$ is open in $X$ for each $y\in Y.$ The
following lemma will be crucial in the proofs.

\begin{lemma}
(Yannelis and Prabhakar, \cite{yan}). \textit{Let }$X$\textit{\ be a
paracompact Hausdorff topological space and }$Y$\textit{\ be a Hausdorff
topological vector space. Let }$T:X\rightarrow 2^{Y}$\textit{\ be a
correspondence with nonempty convex values\ and for each }$y\in Y$\textit{, }%
$T^{-1}(y)$\textit{\ is open in }$X$\textit{. Then, }$T$\textit{\ has a
continuous selection that is, there exists a continuous function }$%
f:X\rightarrow Y$\textit{\ such that }$f(x)\in T(x)$\textit{\ for each }$%
x\in X$\textit{.\medskip }
\end{lemma}

$T$ is said to be \textit{upper semicontinuous,} if for each $x\in X$ and
each open set $V$ in $Y$ with $T(x)\subset V$, there exists an open
neighborhood $U$ of $x$ in $X$ such that $T(x)\subset V$ for each $y\in U$. $%
T$ is said to be \textit{lower semicontinuous,} if for each x$\in X$ and
each open set $V$ in $Y$ with $T(x)\cap V\neq \emptyset $, there exists an
open neighborhood $U$ of $x$ in $X$ such that $T(y)\cap V\neq \emptyset $
for each $y\in U$. $T$ is said to have \textit{open lower sections} if $%
T^{-1}(y)$ is open in $X$ for each $y\in Y.\medskip $

We recall the following definitions concerning the generalized convexity of
the correspondences. The reader is referred to \cite{lin2007}.

Let $X$ be a nonempty convex subset of a topological vector space $E_{i}$, $%
Y $ be a topological vector space and let $C:X\rightarrow 2^{Y}$ be a
correspondence such that for each $x\in X,$ $C(x)$ is a proper, closed and
convex cone with int$C(x)\neq \emptyset .$

A correspondence $T:X\rightarrow 2^{Y}$ is said to be \textit{concave} if
for any $x_{1},x_{2}\in $ $X$ and $\lambda \in \lbrack 0,1],$ $\lambda
T(x_{1})+(1-\lambda )T(x_{2})\subseteq T(\lambda x_{1}+(1-\lambda )x_{2}).$

Let $x\in X$ be arbitrary and let $T:X\times X\rightarrow 2^{Y}$ be a
correspondence.

a) $T$ is said to be $C(x)-$\textit{quasi-convex} if for any $y_{1},y_{2}\in 
$ $X$ and $\lambda \in \lbrack 0,1],$ $\ $\ we have either

$T(x,y_{1})\subseteq T(x,\lambda y_{1}+(1-\lambda )y_{2})+C(x)$ or $%
T(x,y_{2})\subseteq T(x,\lambda y_{1}+(1-\lambda )y_{2})+C(x);$

a) $T$ is said to be $C(x)-$\textit{quasi-convex-like} if for any $%
y_{1},y_{2}\in $ $X$ and $\lambda \in \lbrack 0,1],$ $\ $\ we have either

$T(x,\lambda y_{1}+(1-\lambda )y_{2})\subseteq T(x,y_{1})-C(x)$ or $%
T(x,\lambda y_{1}+(1-\lambda )y_{2})\subseteq T(x,y_{2})-C(x);$

c) $T$ is said to be natural $C(x)-$\textit{quasi-concave} if for any $%
y_{1},y_{2}\in $ $X$ , $\lambda \in \lbrack 0,1]$ and $z_{1}\in T(x,y_{1}),$ 
$z_{2}\in T(x,y_{2}),$ there exists $z\in T(x,\lambda y_{1}+(1-\lambda
)y_{2})$ such that $z\in $co$\{z_{1},z_{2}\}-C(x);$

d) $T$ is said to be $C(x)-$\textit{convex on }$X$ if for any $%
y_{1},y_{2}\in $ $X$ and $\lambda \in \lbrack 0,1],$

$T(x,\lambda y_{1}+(1-\lambda )y_{2})\subseteq \lambda T(x,y_{1})+(1-\lambda
)T(x,y_{2})-C(x).$\medskip

The rest of the paper is organized as follows: Section 2 contains existence
results for the equilibrium of generalized games in choice form and of
generalized abstract economies. The existence of solutions for systems of
vector quasi-equilibrium problems is studied in Section 3. The conclusions
are added in Section 4. All proofs are collected in an Appendix.

\section{Equilibrium results}

In this section, we will consider the family $(X_{i})_{i\in N}$ of
individual strategies and the families $(A_{i})_{i\in N}$ $\ $and $%
(F_{i})_{i\in N}$ of constraint correspondences, where $A_{i}:X_{-i}%
\rightarrow 2^{X_{i}}$ and $F_{i}:X_{-i}\rightarrow 2^{X_{i}}$ for each $%
i\in N$. Let $X=\prod\nolimits_{i\in I}X_{i}$.

We\textit{\ }denote $x_{-i}=(x_{1},...,x_{i-1},x_{i+1},...,x_{n})$, $%
X_{-i}=\prod\nolimits_{i\neq j}X_{i}$ and $(x_{-i},X_{i})=%
\{(x_{-i},x_{i}):x_{i}\in X_{i}\}.\medskip $

We introduce the following definitions, which generalize the ones due to
Ferrara and Stefanescu (\cite{fe}).\medskip\ These authors mentioned, in
order to motivate the presentation of their new model, that the notion of
equilibrium in choice does not require the players' preferences to be
explicitly represented, or the choice of a player to be the best reply to
the strategy combination of the others. As in the model due to Kim and Tan 
\cite{kim}, the correspondences $(F_{i})_{i\in N}$ represent fuzzy
constraints.

The existence of two constraint correspondences is perfectly justified by
experience. For instance, in a recent paper (\cite{he}), Herv\'{e}s-Beloso
and Patriche provided relevant examples of abstract economies with two
constraint correspondences. The model of Kim and Tan \cite{kim} also
contains an additional "fuzzy" constraint.\medskip

\begin{definition}
We define a \textit{choice profile} \textit{under restrictions} as any
collection $\mathcal{C}:=(\mathcal{C}_{i})_{i\in N}$ of nonempty subsets of $%
X\times X$ such that $\mathcal{C}_{i}\subset $Gr$A_{i}\times $Gr$F_{i}$ for
each $i\in N.$
\end{definition}

\begin{definition}
A \textit{generalized} \textit{game in choice form} is a family $%
((X_{i})_{i\in N},$ $(A_{i})_{i\in N},$ $(F_{i})_{i\in N},$ $(\mathcal{C}%
_{i})_{i\in N})$, where $\mathcal{C}=(\mathcal{C}_{i})_{i\in N\text{ }}$is a
choice profile under restrictions.
\end{definition}

\begin{definition}
An \textit{equilibrium in choice} of a generalized game in choice form $%
((X_{i})_{i\in N},(A_{i})_{i\in N}\ ,(F_{i})_{i\in N},$ $(\mathcal{C}%
_{i})_{i\in N})$ is any game strategy $\ (x^{\ast },y^{\ast })$ with the
property
\end{definition}

$\forall i\in N,$ $((x_{-i}^{\ast },A_{i}(x_{-i}^{\ast })),(y_{-i}^{\ast
},F_{i}(y_{-i}^{\ast }))\cap \mathcal{C}_{i}\neq \emptyset \Rightarrow
(x^{\ast },y^{\ast })\in \mathcal{C}_{i}.$

If $(x^{\ast },y^{\ast })\in \cap _{i\in N}\mathcal{C}_{i},$ then $(x^{\ast
},y^{\ast })$ is said to be a \textit{strong equilibrium in choice.\medskip }

We denote $\mathcal{C}(x_{-i},y_{-i})$ the upper section through $%
(x_{-i},y_{-i})$ of a set $\mathcal{C}\subset X\times X,$ i.e., $\mathcal{C}%
(x_{-i},y_{-i})=\{(x_{i},y_{i})\in X_{i}\times
X_{i}:((x_{-i},x_{i}),(y_{-i},y_{i}))\in \mathcal{C}\}$ and $\mathcal{C}%
(x_{i},y_{i})$ the lower section through $(x_{i},y_{i})$ of the set $%
\mathcal{C},$ i.e., $\mathcal{C}(x_{i},y_{i})=\{(x_{-i},y_{-i})\in
X_{-i}\times X_{-i}:((x_{-i},x_{i}),(y_{-i},y_{i}))\in \mathcal{C}\}$ $.$

We will make the following assumption:\newline
(A) for each $(x,y)\in X\times X,$ there exists $i\in N$ such that $\mathcal{%
C}_{i}(x_{-i},y_{-i})\neq \emptyset .\medskip $

Finally, we note that we determine the existence of strong equilibrium in
choice for all the situations considered in this section, if we suppose, in
addition, that $\mathcal{C}_{i}(x_{-i})\neq \emptyset $ for each $x_{-i}\in
X_{-i}.\medskip $

Further, we study the existence of equilibrium for a generalized game in
choice form. All proofs will be collected in an Appendix, at the end of the
paper.\medskip\ We underline that our approach differs very much from the
one due to Ferrara and Stefanescu in \cite{fe} or to Stefanescu, Ferrara and
Stefanescu in \cite{st}.\medskip

This is our first result. The assumptions refer mainly to the topological
properties of the upper and lower sections of the sets $\mathcal{C}_{i},$ $%
i\in N.$

\begin{theorem}
Let $(X_{i},A_{i},F_{i},\mathcal{C}_{i})_{i\in N}$ be a generalized \textit{%
game in choice form. Assume that, for each }$i\in N,$ the following
conditions are fulfilled:
\end{theorem}

\textit{a)} $X_{i}$ \textit{is a nonempty, convex and compact set in a} 
\textit{Hausdorff locally convex space }$E_{i}$;

b) $W_{i}=\{(x_{-i},y_{-i})\in X_{-i}\times X_{-i}:A_{i}(x_{-i})\times
F_{i}(y_{-i})\cap \mathcal{C}_{i}(x_{-i},y_{-i})\neq \emptyset \}\mathit{\ }$
\textit{is nonempty and closed}$;$

\textit{c) }$\mathcal{C}_{i}(x_{i},y_{i})$ \textit{is open for} \textit{each 
} $(x_{i},y_{i})\in X_{i}\times X_{i};$

\textit{d)} $\mathcal{C}_{i}(x_{-i},y_{-i})\mathit{\ }$ \textit{is convex or
empty for each} $(x_{-i},y_{-i})\in X_{-i}\times X_{-i}.$

\textit{Then, the game admits equilibria in choice.\medskip }

As a consequence, we obtain the following result concerning the existence of
the strong equilibrium in choice.

\begin{theorem}
Let $(X_{i},A_{i},F_{i},\mathcal{C}_{i})_{i\in N}$ be a generalized \textit{%
game in choice form. Assume that, for each }$i\in N,$ the following
conditions are fulfilled:
\end{theorem}

\textit{a)} $X_{i}$ \textit{is a nonempty, convex and compact set in a} 
\textit{Hausdorff locally convex space }$E_{i}$;

\textit{b) }$(A_{i}(x_{-i})\times F_{i}(y_{-i}))\cap \mathcal{C}%
_{i}(x_{-i},y_{-i})\mathit{\ }$ \textit{is nonempty for each} $%
(x_{-i},y_{-i})\in X_{-i}\times X_{-i}.$

\textit{c) }$\mathcal{C}_{i}(x_{i},y_{i})$ \textit{is open for} \textit{each 
} $(x_{i},y_{i})\in X_{i}\times X_{i};$

\textit{d)} $\mathcal{C}_{i}(x_{-i},y_{-i})\mathit{\ }$ \textit{is convex
for each} $x_{-i}$ \textit{and} $y_{-i}\in X_{-i}.$

\textit{Then, the generalized game admits strong equilibria in
choice.\medskip }

Further, we recall that the model of a generalized abstract economy, which
was introduced by Kim and Tan in \cite{kim}, is more general then the one
due to Shafer and Sonnenschein \cite{sh}. In Kim and Tan's model, each agent 
$i$ has, in addition, a fuzzy constraint correspondence $F_{i}$.

In order to prove the existence of equilibrium for the generalized abstract
economies, we attach an auxiliary generalized game in choice form and we use
the above theorems. The strong equilibrium in choice for the generalized
games in choice form will be equilibrium pairs for generalized abstract
economies.\medskip

Let $(X_{i})_{i\in N}$ be the family of the individual sets of strategies
and let $X=\prod\nolimits_{i\in I}X_{i}$.\medskip

We will work with the following variant of a \textit{generalized abstract
economy:\ }

$\Gamma =(X_{i},A_{i},F_{i},P_{i})_{i\in N}$ is defined as a family of
quadruplets $(X_{i},A_{i},F_{i},P_{i}),$ where for each $i\in N$, $%
P_{i}:X\times X\rightarrow 2^{X_{i}}$ is a preference\ correspondence and $%
A_{i},F_{i}:X_{-i}\rightarrow 2^{X_{i}}$ are constraint\ correspondences.

An \textit{equilibrium} for $\Gamma $ is a pair $(x^{\ast },y^{\ast })\in $ $%
X\times X$ which satisfies for each $i\in N:$ $\ x_{i}^{\ast }\in \overline{A%
}_{i}(x_{-i}^{\ast }),$ $y_{i}^{\ast }\in \overline{F}_{i}(y_{-i}^{\ast })$
and $A_{i}(x_{-i}^{\ast })\cap P_{i}(x^{\ast },y^{\ast })=\emptyset
.\medskip $

Theorem 3 introduces new conditions which ensure the existence of equilibria
for abstract economies. Its proof is based on Theorem 2.

\begin{theorem}
Let $(X_{i},A_{i},F_{i},P_{i})_{i\in N}$ be a generalized abstract economy%
\textit{. Assume that, for each }$i\in N,$ the following conditions are
fulfilled:
\end{theorem}

\textit{a)} $X_{i}$ \textit{is a nonempty, convex and compact set in a} 
\textit{Hausdorff locally convex space }$E_{i}$;

\textit{b)} $A_{i},F_{i}$\textit{\ have nonempty, convex values and open
lower sections;}

\textit{c)} \textit{the set }$\{(x_{i},y_{i})\in X_{i}\times X_{i}:$ $%
A_{i}(x_{-i})\cap P_{i}(x,y)=\emptyset \}\cap (A_{i}(x_{-i})\times
F_{i}(y_{-i}))$ \ \textit{is nonempty for each} $(x_{-i},y_{-i})\in
X_{-i}\times X_{-i}$\textit{;}

\textit{d) }$\{(x_{i},y_{i})\in X_{i}\times X_{i}:$\textit{\ }$%
A_{i}(x_{-i})\cap P_{i}(x,y)=\emptyset \}$\textit{\ is convex for each }$%
(x_{-i},y_{-i})\in X_{-i}\times X_{-i};$

\textit{e)} $\{(x_{-i},y_{-i})\in X_{-i}\times X_{-i}:A_{i}(x_{-i})\cap
P_{i}(x,y)=\emptyset \}\mathit{\ }$ \textit{is open for each} $%
(x_{i},y_{i})\in X_{i}\times X_{i}.$

\textit{Then, there exists a pair }$(x^{\ast },y^{\ast })\in $\textit{\ }$%
X\times X$\textit{\ such that }$\ x_{i}^{\ast }\in A_{i}(x_{-i}^{\ast }),$%
\textit{\ }$y_{i}^{\ast }\in F_{i}(y_{-i}^{\ast })$\textit{\ and }$%
A_{i}(x_{-i}^{\ast })\cap P_{i}(x^{\ast },y^{\ast })=\emptyset $\textit{\
for each} $i\in N.\medskip $

We obtain the following result as a direct consequence of the above theorem.
It will be used in the next section in order to prove the existence of
solutions for systems of vector quasi-equilibrium problems.

\begin{theorem}
Let $(X_{i},A_{i},F_{i},P_{i})_{i\in N}$ be a generalized abstract economy%
\textit{. Assume that, for each }$i\in N,$ the following conditions are
fulfilled:
\end{theorem}

\textit{a)} $X_{i}$ \textit{is a nonempty, convex and compact set in a} 
\textit{Hausdorff locally convex space }$E_{i}$;

\textit{b)} $A_{i},F_{i}$\textit{\ have nonempty, convex values and open
lower sections;}

\textit{c)} \textit{the set }$\{(x_{i},y_{i})\in X_{i}\times X_{i}:$ $%
\overline{A}_{i}(x_{-i})\cap P_{i}(x,y)=\emptyset \}\cap
(A_{i}(x_{-i})\times F_{i}(y_{-i}))$ \ \textit{is nonempty for each} $%
(x_{-i},y_{-i})\in X_{-i}\times X_{-i}$\textit{;}

\textit{d) }$\{(x_{i},y_{i})\in X_{i}\times X_{i}:$\textit{\ }$\overline{A}%
_{i}(x_{-i})\cap P_{i}(x,y)=\emptyset \}$\textit{\ is convex for each }$%
(x_{-i},y_{-i})\in X_{-i}\times X_{-i};$

\textit{e)} $\{(x_{-i},y_{-i})\in X_{-i}\times X_{-i}:\overline{A}%
_{i}(x_{-i})\cap P_{i}(x,y)=\emptyset \}\mathit{\ }$ \textit{is open for each%
} $(x_{i},y_{i})\in X_{i}\times X_{i}.$

\textit{Then, there exists a pair }$(x^{\ast },y^{\ast })\in $\textit{\ }$%
X\times X$\textit{\ such that }$\ x_{i}^{\ast }\in A_{i}(x_{-i}^{\ast }),$%
\textit{\ }$y_{i}^{\ast }\in F_{i}(y_{-i}^{\ast })$\textit{\ and }$%
A_{i}(x_{-i}^{\ast })\cap P_{i}(x^{\ast },y^{\ast })=\emptyset $\textit{\
for each} $i\in N.\medskip $

We note that, according to the above theorem, there exists a pair $(x^{\ast
},y^{\ast })\in $ $X\times X$ such that $\ x_{i}^{\ast }\in
A_{i}(x_{-i}^{\ast }),$ $y_{i}^{\ast }\in F_{i}(y_{-i}^{\ast })$ and $%
\overline{A}_{i}(x_{-i}^{\ast })\cap P_{i}(x^{\ast },y^{\ast })=\emptyset $
for each $i\in N.$ Obviously, $\overline{A}_{i}(x_{-i}^{\ast })\cap
P_{i}(x^{\ast },y^{\ast })=\emptyset $ implies $A_{i}(x_{-i}^{\ast })\cap
P_{i}(x^{\ast },y^{\ast })=\emptyset .\medskip $

Theorem 5 states the existence of equilibrium for a generalized abstract
economy, in which the constraint correspondences are lower semicontinuous
and the preference correspondences have open graphs. No continuity
assumptions are made over the fuzzy constraint correspondences. No
assumptions are made over the values of the preference correspondences. The
values of the constraint correspondences needn't fulfill topological
conditions. The proof is based on the construction, for each player, of a
new type of correspondence, which satisfies the Kakutani fixed point Theorem.

\begin{theorem}
Let $(X_{i},A_{i},F_{i},P_{i})_{i\in N}$ be an abstract economy\textit{.
Assume that, for each }$i\in N,$ the following conditions are fulfilled:
\end{theorem}

\textit{a)} $X_{i}$ \textit{is a nonempty, convex and compact set in a} 
\textit{Hausdorff locally convex space }$E_{i}$;

\textit{b) }$A_{i}$\textit{\ is lower semicontinuous; {}}

\textit{c) }$\overline{A_{i}}$\textit{\ and }$\overline{F_{i}}$ \textit{are
nonempty and convex valued;}

\textit{d) }$P_{i}$\textit{\ has an open graph;}

\textit{e)} $\{(x_{i},y_{i})\in X_{i}\times X_{i}:A_{i}(x_{-i})\cap
P_{i}(x,y)=\emptyset \}\mathit{\ }$ \textit{is convex and nonempty for each} 
$(x_{-i},y_{-i})\in X_{-i}\times X_{-i}.$

\textit{Then, there exists a pair }$(x^{\ast },y^{\ast })\in $\textit{\ }$%
X\times X$\textit{\ which satisfies for each }$i\in N:$\textit{\ }$\
x_{i}^{\ast }\in \overline{A_{i}}(x_{-i}^{\ast }),$ $y_{i}^{\ast }\in 
\overline{F_{i}}(y_{-i}^{\ast })$\textit{\ and} $A_{i}(x_{-i}^{\ast })\cap
P_{i}(x^{\ast },y^{\ast })=\emptyset .$

\begin{remark}
\textit{Assumption d) of the above theorem implies that if} $(x,y)\in
X\times X,$ $x_{i}\notin P_{i}(x,y).\medskip $
\end{remark}

\section{Systems of vector quasi-equilibrium problems}

We start this section with the presentation of the problem we
approach.\medskip

For each $i\in N,$ let $X_{i}$ be a nonempty subset of a Hausdorff
topological vector space $E_{i},$ $Y_{i}$ a topological vector space and let 
$X=\tprod\limits_{i\in N}X_{i}$ $.$ For each $i\in N,$ let $A_{i},$ $%
F_{i}:X_{-i}\rightarrow 2^{X_{i}}$, $C_{i}:X_{-i}\rightarrow 2^{Y_{i}}$ and $%
f_{i}:X\times X\times X_{i}\rightarrow 2^{Y_{i}}$ be correspondences with
nonempty values. We consider the following systems of generalized vector
quasi-equilibrium problems (in short, SGVQEP):

SGVQEP (I) Find $(x^{\ast },y^{\ast })\in X\times X$ such that for each $%
i\in N,$ $x_{i}^{\ast }\in \overline{A}_{i}(x_{-i}^{\ast }),$ $y_{i}^{\ast
}\in \overline{F}_{i}(y_{-i}^{\ast })$ and $f_{i}(x^{\ast },y^{\ast
},u_{i})\subseteq C_{i}(x_{-i}^{\ast })$ for each $u_{i}\in
A_{i}(x_{-i}^{\ast }).$

SGVQEP (II) Find $(x^{\ast },y^{\ast })\in X\times X$ such that for each $%
i\in N,$ $x_{i}^{\ast }\in \overline{A}_{i}(x_{-i}^{\ast }),$ $y_{i}^{\ast
}\in \overline{F}_{i}(y_{-i}^{\ast })$ and $f_{i}(x^{\ast },y^{\ast
},u_{i})\cap C_{i}(x_{-i}^{\ast })\neq \emptyset $ for each $u_{i}\in
A_{i}(x_{-i}^{\ast }).$

SGVQEP (III) Find $(x^{\ast },y^{\ast })\in X\times X$ such that for each $%
i\in N,$ $x_{i}^{\ast }\in \overline{A}_{i}(x_{-i}^{\ast }),$ $y_{i}^{\ast
}\in \overline{F}_{i}(y_{-i}^{\ast })$ and $f_{i}(x^{\ast },y^{\ast
},u_{i})\cap (-$int$C_{i}(x_{-i}^{\ast }))=\emptyset $ for each $u_{i}\in
A_{i}(x_{-i}^{\ast }).$

SGVQEP (IV) Find $(x^{\ast },y^{\ast })\in X\times X$ such that for each $%
i\in N,$ $x_{i}^{\ast }\in \overline{A}_{i}(x_{-i}^{\ast }),$ $y_{i}^{\ast
}\in \overline{F}_{i}(y_{-i}^{\ast })$ and $f_{i}(x^{\ast },y^{\ast
},u_{i})\nsubseteq $int$C_{i}(x_{-i}^{\ast })$ for each $u_{i}\in
A_{i}(x_{-i}^{\ast }).\medskip $

Further, we will work in the following setting:

For each $i\in N,$ let $X_{i}$ be a nonempty subset of a Hausdorff locally
convex space $E_{i},$ $Y_{i}$ a topological vector space and let $%
X=\tprod\limits_{i\in N}X_{i}$ $.$ For each $i\in N,$ let $A_{i},$ $%
F_{i}:X_{-i}\rightarrow 2^{X_{i}}$, $C_{i}:X_{-i}\rightarrow 2^{Y_{i}}$ and $%
f_{i}:X\times X\times X_{i}\rightarrow 2^{Y_{i}}$ be correspondences with
nonempty values. Suppose that for each $x_{-i}\in X_{-i},$ $C_{i}(x_{-i})$
is a proper, closed and convex cone with int$C_{i}(x_{-i})\neq \emptyset
.\medskip $

We obtain the following results which improve the theorems from \cite%
{lin2007}. All proofs are included in Appendix. We mention that the proofs
are different from the ones given in \cite{lin2007}.\medskip

Now, we present two different theorems, which concern the existence of
solutions for SGVQEP (I). Both of them are established under hypotheses
which are weaker than those from \cite{lin2007}.

We note that, in Theorem 6, the upper semicontinuity of each $f_{i}$ was
restricted to the upper semicontinuity of each \textit{\ }$f_{i}(\cdot
,x_{i},\cdot ,y_{i},\cdot ):$\textit{\ }$X_{-i}\times X_{-i}\times
X_{i}\rightarrow 2^{Y_{i}}.$ The convexity of the values of the
correspondences $f_{i}$ is not assumed.

\begin{theorem}
\textit{For each }$i\in N,$\textit{\ let }$f_{i}:X\times X\times
X_{i}\rightarrow 2^{Y_{i}}$\textit{\ be a correspondence with nonempty and
closed values. For each }$i\in N,$\textit{\ suppose that:}
\end{theorem}

\textit{a) }$A_{i},F_{i}$\textit{\ have nonempty, convex values and open
lower sections;}

\textit{b) }$f_{i}(\cdot ,x_{i},\cdot ,y_{i},\cdot ):$\textit{\ }$%
X_{-i}\times X_{-i}\times X_{i}\rightarrow 2^{Y_{i}}$\textit{\ is upper
semicontinuous for each }$(x_{i},y_{i})\in X_{i}\times X_{i};$

\textit{c) the correspondence }$W_{i}$ $:X_{-i}\rightarrow 2^{Y_{i}}$, 
\textit{defined by }$W_{i}(x_{-i})=Y_{i}\backslash (-$\textit{int}$%
C_{i}(x_{-i}))$\textit{\ for each }$x_{-i}\in X_{-i},$\textit{\ is upper
semicontinuous;}

\textit{d) for each }$(x_{-i},y_{-i})\in X_{-i}\times X_{-i},$ \textit{there
exists }$(x_{i},y_{i})\in A_{i}(x_{-i})\times F_{i}(y_{-i})$\textit{\ such
that }$f_{i}(x,y,\overline{A}_{i}(x_{-i}))\subseteq $int$C_{i}(x_{-i})$%
\textit{\ }$;$

\textit{e) for each }$(x_{-i},y_{-i},u_{i})\in X_{-i}\times X_{-i}\times
X_{i},$\textit{\ }$f_{i}(x_{-i},\cdot ,y_{-i},\cdot ,u_{i}):X_{i}\times
X_{i}\rightarrow 2^{Y_{i}}$\textit{\ is }$(-$int$C_{i}(x_{-i}))-$\textit{\
quasi-convex-like.}

\textit{Then, there exists a pair }$(x^{\ast },y^{\ast })\in X\times X$%
\textit{\ such that for each }$i\in N,$\textit{\ }$x_{i}^{\ast }\in
A_{i}(x_{-i}^{\ast }),$\textit{\ }$y_{i}^{\ast }\in F_{i}(y_{-i}^{\ast })$%
\textit{\ and }$f_{i}(x^{\ast },y^{\ast },u_{i})\subseteq C_{i}(x_{-i}^{\ast
})$\textit{\ for each }$u_{i}\in A_{i}(x_{-i}^{\ast }),$\textit{\ that is, }$%
(x^{\ast },y^{\ast })$\textit{\ is a solution for SGVQEP (I).\medskip }

We will show that the existence of solutions for SGVQEP (I) can be derived
from Theorem 5, following the method used by Lin, Chen and Ansari in \cite%
{lin2007}. Theorem 7 and Theorem 3.2.1 from \cite{lin2007} are comparable.
After the proof, we will present a comparison between our result and the
quoted one.

\begin{theorem}
\textit{For each }$i\in N,$\textit{\ let }$f_{i}:X\times X\times
X_{i}\rightarrow 2^{Y_{i}}$\textit{\ be a lower semicontinuous
correspondence with nonempty values. For each }$i\in N,$\textit{\ suppose
that:}
\end{theorem}

\textit{a) }$A_{i}$\textit{\ is lower semicontinuous;}

\textit{b) }$\overline{A_{i}}$\textit{\ and }$\overline{F_{i}}$\textit{\ are}
\textit{nonempty and convex valued;}

\textit{c) }$C_{i}$ \textit{is an upper semicontinuous correspondence;}

\textit{d) for each }$(x_{-i},y_{-i})\in X_{-i}\times X_{-i},$ $%
\{(x_{i},y_{i})\in X_{i}\times X_{i}:$\textit{\ }$f_{i}(x,y,u_{i})\subseteq
C_{i}(x_{-i})$\textit{\ for each }$u_{i}\in A_{i}(x_{-i})\}$\textit{\ is
nonempty}$;$

\textit{e) for each }$(x_{-i},y_{-i},u_{i})\in X_{-i}\times X_{-i}\times
X_{i},$\textit{\ }$f_{i}(x_{-i},\cdot ,y_{-i},\cdot ,u_{i}):X_{i}\times
X_{i}\rightarrow 2^{Y_{i}}$\textit{\ is }$(-C_{i}(x_{-i}))-$\textit{\
quasi-convex-like.}

\textit{Then, there exists a pair }$(x^{\ast },y^{\ast })\in X\times X$%
\textit{\ such that for each }$i\in N,$\textit{\ }$x_{i}^{\ast }\in 
\overline{A}_{i}(x_{-i}^{\ast }),$\textit{\ }$y_{i}^{\ast }\in \overline{F}%
_{i}(y_{-i}^{\ast })$\textit{\ and }$f_{i}(x^{\ast },y^{\ast
},u_{i})\subseteq C_{i}(x_{-i}^{\ast })$\textit{\ for each }$u_{i}\in
A_{i}(x_{-i}^{\ast }),$\textit{\ that is, }$(x^{\ast },y^{\ast })$\textit{\
is a solution for SGVQEP (I).\medskip }

\begin{remark}
Assumption d) implies that for each $x,y\in X,$ \textit{\ }$%
f_{i}(x,y,x_{i})\subseteq C_{i}(x_{-i}).$
\end{remark}

\begin{remark}
We note that no continuity assumptions was made over the correspondences $%
F_{i},$ $i\in N.$ The functions $f_{i}:X\times X\times X_{i}\rightarrow
2^{Y_{i}}$ are not quasi-concave in the last argument.\medskip
\end{remark}

SGVQEP (II) has solutions if we suppose that the following assumptions hold.

\begin{theorem}
\textit{For each }$i\in N,$\textit{\ let }$f_{i}:X\times X\times
X_{i}\rightarrow 2^{Y_{i}}$\textit{\ be a correspondence with nonempty
values. For each }$i\in N,$\textit{\ suppose that:}
\end{theorem}

\textit{a) }$A_{i},F_{i}$\textit{\ have nonempty, convex values and open
lower sections;}

\textit{b) }$f_{i}(\cdot ,x_{i},\cdot ,y_{i},\cdot ):$\textit{\ }$%
X_{-i}\times X_{-i}\times X_{i}\rightarrow 2^{Y_{i}}$\textit{\ is lower
semicontinuous for each }$(x_{i},y_{i})\in X_{i}\times X_{i};$

\textit{c) the correspondence }$W_{i}$ $:X_{-i}\rightarrow 2^{Y_{i}}$\textit{%
, defined by }$W_{i}(x_{-i})=Y_{i}\backslash (-$\textit{int}$C_{i}(x_{-i}))$%
\textit{\ for each }$x_{-i}\in X_{-i},$\textit{\ is upper semicontinuous;}

\textit{d) there exists }$(x_{i},y_{i})\in A_{i}(x_{-i})\times F_{i}(y_{-i})$%
\textit{\ such that }$f_{i}(x,y,u_{i})\cap $\textit{int}$C_{i}(x_{-i})\neq
\emptyset $\textit{\ \textit{for each} }$u_{i}\in \overline{A}_{i}(x_{-i})$ 
\textit{and} \textit{for each} $(x_{-i},y_{-i})\in X_{-i}\times X_{-i};$

\textit{e) for each }$(x_{-i},y_{-i},u_{i})\in X_{-i}\times X_{-i}\times
X_{i},$\textit{\ }$f_{i}(x_{-i},\cdot ,y_{-i},\cdot ,u_{i}):X_{i}\times
X_{i}\rightarrow 2^{Y_{i}}$\textit{\ is natural int}$C_{i}(x_{-i})-$\textit{%
\ quasi-concave.}

\textit{Then, there exists a pair }$(x^{\ast },y^{\ast })\in X\times X$%
\textit{\ such that for each }$i\in N,$\textit{\ }$x_{i}^{\ast }\in
A_{i}(x_{-i}^{\ast }),$\textit{\ }$y_{i}^{\ast }\in F_{i}(y_{-i}^{\ast })$%
\textit{\ and } $f_{i}(x^{\ast },y^{\ast },u_{i})\cap C_{i}(x_{-i}^{\ast
})\neq \emptyset $\textit{\ for each }$u_{i}\in A_{i}(x_{-i}^{\ast }),$%
\textit{\ that is, }$(x^{\ast },y^{\ast })$\textit{\ is a solution for
SGVQEP (II).\medskip }

The following theorem gives sufficient conditions for the solvability of
SGVQEP (II), given upper semicontinuous correspondences $f_{i}:X\times
X\times X_{i}\rightarrow 2^{Y_{i}}$, $i\in N.$

\begin{theorem}
For each $i\in N,$ let $f_{i}:X\times X\times X_{i}\rightarrow 2^{Y_{i}}$ be
an upper semicontinuous correspondence with nonempty values. For each $i\in
N,$ suppose that:
\end{theorem}

\textit{a) }$A_{i}$\textit{\ is lower semicontinuous;}

\textit{b) }$\overline{A_{i}}$\textit{\ and }$\overline{F_{i}}$\textit{\ are}
\textit{nonempty and convex valued;}

\textit{c) }$C_{i}$ \textit{is an upper semicontinuous correspondence;}

\textit{d) for each}\ $(x_{-i},y_{-i})\in X_{-i}\times X_{-i},$ $%
\{(x_{i},y_{i})\in X_{i}\times X_{i}:$ $f_{i}(x,y,u_{i})\cap
C_{i}(x_{-i})\neq \emptyset $ \textit{for each }$u_{i}\in A_{i}(x_{-i})\}$ 
\textit{is nonempty}$;$

\textit{e) for each} $(x_{-i},y_{-i},u_{i})\in X_{-i}\times X_{-i}\times
X_{i},$ $f_{i}(x_{-i},\cdot ,y_{-i},\cdot ,u_{i}):X_{i}\times
X_{i}\rightarrow 2^{Y_{i}}$ \textit{is natural }$C_{i}(x_{-i})-$\textit{\
quasi-concave.}

\textit{Then, there exists a pair }$(x^{\ast },y^{\ast })\in X\times X$%
\textit{\ such that for each }$i\in N,$\textit{\ }$x_{i}^{\ast }\in 
\overline{A}_{i}(x_{-i}^{\ast }),$\textit{\ }$y_{i}^{\ast }\in \overline{F}%
_{i}(y_{-i}^{\ast })$\textit{\ and } $f_{i}(x^{\ast },y^{\ast },u_{i})\cap
C_{i}(x_{-i}^{\ast })\neq \emptyset $\textit{\ for each }$u_{i}\in
A_{i}(x_{-i}^{\ast }),$\textit{\ that is, }$(x^{\ast },y^{\ast })$\textit{\
is a solution for SGVQEP (II).\medskip }

The existence of the solutions for SGVQEP (III)\textit{\ }\ is stated below.

\begin{theorem}
\textit{For each }$i\in N,$\textit{\ let }$f_{i}:X\times X\times
X_{i}\rightarrow 2^{Y_{i}}$\textit{\ be a correspondence with nonempty,
closed values. For each }$i\in N,$\textit{\ suppose that:}
\end{theorem}

\textit{a) }$A_{i},F_{i}$\textit{\ have nonempty, convex values and open
lower sections;}

\textit{b) }$f_{i}(\cdot ,x_{i},\cdot ,y_{i},\cdot ):$\textit{\ }$%
X_{-i}\times X_{-i}\times X_{i}\rightarrow 2^{Y_{i}}$\textit{\ is upper
semicontinuous for each }$(x_{i},y_{i})\in X_{i}\times X_{i};$

\textit{c) }$C_{i}$ \textit{is an upper semicontinuous correspondence;}

\textit{d) there exists }$(x_{i},y_{i})\in A_{i}(x_{-i})\times F_{i}(y_{-i})$%
\textit{\ such that }$f_{i}(x,y,u_{i})\cap (-C_{i}(x_{-i}))=\emptyset $%
\textit{\ for each }$\mathit{u}_{i}\in \overline{A}_{i}(x_{-i})$ \textit{and
for each} $(x_{-i},y_{-i})\in X_{-i}\times X_{-i};$

\textit{e) for each }$(x_{-i},y_{-i},u_{i})\in X_{-i}\times X_{-i}\times
X_{i},$\textit{\ }$f_{i}(x_{-i},\cdot ,y_{-i},\cdot ,u_{i}):X_{i}\times
X_{i}\rightarrow 2^{Y_{i}}$\textit{\ is }$(-C_{i}(x_{-i}))-$\ \textit{%
quasi-convex-like.}

\textit{Then, there exists a pair }$(x^{\ast },y^{\ast })\in X\times X$%
\textit{\ such that for each }$i\in N,$\textit{\ }$x_{i}^{\ast }\in
A_{i}(x_{-i}^{\ast }),$\textit{\ }$y_{i}^{\ast }\in F_{i}(y_{-i}^{\ast })$%
\textit{\ and } $f_{i}(x^{\ast },y^{\ast },u_{i})\cap (-$int$%
C_{i}(x_{-i}^{\ast }))=\emptyset $\textit{\ for each }$u_{i}\in
A_{i}(x_{-i}^{\ast }),$\textit{\ that is, }$(x^{\ast },y^{\ast })$\textit{\
is a solution for SGVQEP (III).\medskip }

The question, whether Theorem 5 implies that SGVQEP (III) has solutions,%
\textit{\ }arises naturally. Assumptions on $f_{i}:X\times X\times
X_{i}\rightarrow 2^{Y_{i}}$ refers to lower semicontinuity.

\begin{theorem}
For each $i\in N,$ let $f_{i}:X\times X\times X_{i}\rightarrow 2^{Y_{i}}$ be
a lower semicontinuous correspondence with nonempty values. For each $i\in
N, $ suppose that:
\end{theorem}

\textit{a) }$A_{i}$\textit{\ is lower semicontinuous;}

\textit{b) }$\overline{A_{i}}$\textit{\ and }$\overline{F_{i}}$ \textit{are
nonempty and convex valued;}

\textit{c) the correspondence }$W_{i}$ $:X_{-i}\rightarrow 2^{Y_{i}}$, 
\textit{defined by }$W_{i}(x_{-i})=Y_{i}\backslash (-$\textit{int}$%
C_{i}(x_{-i}))$\textit{\ for each }$x_{-i}\in X_{-i},$\textit{\ is upper
semicontinuous;}

\textit{d) for each }$(x_{-i},y_{-i})\in X_{-i}\times X_{-i},$\textit{\ }$%
\{(x_{i},y_{i})\in X_{i}\times X_{i}:$\textit{\ }$f_{i}(x,y,u_{i})\cap (-$%
\textit{int}$C_{i}(x_{-i}))=\emptyset $\textit{\ for each }$u_{i}\in
A_{i}(x_{-i})\}$\textit{\ is nonempty}$;$

\textit{e) for each }$(x_{-i},y_{-i},u_{i})\in X_{-i}\times X_{-i}\times
X_{i},$\textit{\ }$f_{i}(x_{-i},\cdot ,y_{-i},\cdot ,u_{i}):X_{i}\times
X_{i}\rightarrow 2^{Y_{i}}$\textit{\ is }$(-$\textit{int}$C_{i}(x_{-i}))-$%
\textit{\ quasi-convex-like.}

\textit{Then, there exists a pair }$(x^{\ast },y^{\ast })\in X\times X$%
\textit{\ such that for each }$i\in N,$\textit{\ }$x_{i}^{\ast }\in 
\overline{A}_{i}(x_{-i}^{\ast }),$\textit{\ }$y_{i}^{\ast }\in \overline{F}%
_{i}(y_{-i}^{\ast })$\textit{\ and } $f_{i}(x^{\ast },y^{\ast },u_{i})\cap (-
$int$C_{i}(x_{-i}^{\ast }))=\emptyset $\textit{\ for each }$u_{i}\in
A_{i}(x_{-i}^{\ast }),$\textit{\ that is, }$(x^{\ast },y^{\ast })$\textit{\
is a solution for SGVQEP (III).\medskip }

\begin{remark}
The above theorem can be compared with Theorem 3.2.3 from \cite{lin2007}.
\end{remark}

The existence of solutions for SGVQEP (IV\textit{) }is studied at the end of
the paper.

\begin{theorem}
\textit{For each }$i\in N,$\textit{\ let }$f_{i}:X\times X\times
X_{i}\rightarrow 2^{Y_{i}}$\textit{\ be a correspondence with nonempty
values. For each }$i\in N,$\textit{\ suppose that:}
\end{theorem}

\textit{a) }$A_{i},F_{i}$\textit{\ have convex values and open lower
sections;}

\textit{b) }$f_{i}(\cdot ,x_{i},\cdot ,y_{i},\cdot ):$\textit{\ }$%
X_{-i}\times X_{-i}\times X_{i}\rightarrow 2^{Y_{i}}$\textit{\ is lower
semicontinuous for each }$(x_{i},y_{i})\in X_{i}\times X_{i};$

\textit{c) }$C_{i}$ \textit{is an upper semicontinuous correspondence;}

\textit{e) there exists }$(x_{i},y_{i})\in A_{i}(x_{-i})\times F_{i}(y_{-i})$%
\textit{\ such that }$f_{i}(x,y,u_{i})\nsubseteq C_{i}(x_{-i})$ \textit{for
each }$u_{i}\in \overline{A}_{i}(x_{-i})$ \textit{and for each }$%
(x_{-i},y_{-i})\in X_{-i}\times X_{-i};$

\textit{f) for each }$(x_{-i},y_{-i},u_{i})\in X_{-i}\times X_{-i}\times
X_{i},$\textit{\ }$f_{i}(x_{-i},\cdot ,y_{-i},\cdot ,u_{i}):X_{i}\times
X_{i}\rightarrow 2^{Y_{i}}$\textit{\ is }$C_{i}(x_{-i})-$\ \textit{%
quasi-convex.}

\textit{Then, there exists a pair }$(x^{\ast },y^{\ast })\in X\times X$%
\textit{\ such that for each }$i\in N,$\textit{\ }$x_{i}^{\ast }\in
A_{i}(x_{-i}^{\ast }),$\textit{\ }$y_{i}^{\ast }\in F_{i}(y_{-i}^{\ast })$%
\textit{\ and } $f_{i}(x^{\ast },y^{\ast },u_{i})\nsubseteq $int$%
C_{i}(x_{-i}^{\ast })$\textit{\ for each }$u_{i}\in A_{i}(x_{-i}^{\ast }),$%
\textit{\ that is, }$(x^{\ast },y^{\ast })$\textit{\ is a solution for
SGVQEP (IV).\medskip }

In \cite{lin2007}, Theorem 3.2.4 concerning the existence of solutions of
SGVQEP (IV) has been established. We prove that this statement remains valid
without assuming the lower semicontinuity of the correspondences $F_{i},$ $%
i\in N$ and the quasi-convexity of $f_{i}$ in the last argument. Instead, we
assume that \textit{\ }for each $(x_{-i},y_{-i},u_{i})\in X_{-i}\times
X_{-i}\times X_{i},$\ $f_{i}(x_{-i},\cdot ,y_{-i},\cdot ,u_{i}):X_{i}\times
X_{i}\rightarrow 2^{Y_{i}}$\ is int$C_{i}(x_{-i})-$\ quasi-convex. The proof
is also based on Theorem 5.

\begin{theorem}
\textit{For each }$i\in N,$\textit{\ let }$f_{i}:X\times X\times
X_{i}\rightarrow 2^{Y_{i}}$\textit{\ be an upper semicontinuous
correspondence with nonempty values. For each }$i\in N,$\textit{\ suppose
that:}
\end{theorem}

\textit{a) }$A_{i}$\textit{\ is lower semicontinuous;}

\textit{b) }$\overline{A_{i}}$\textit{\ and }$\overline{F_{i}}$ \textit{are
nonempty and convex valued;}

\textit{c) the correspondence }$W_{i}$\textit{\ }$:X_{-i}\rightarrow
2^{Y_{i}}$\textit{, defined by }$W_{i}(x_{-i})=Y_{i}\backslash $\textit{int}$%
C_{i}(x_{-i})$\textit{\ for each }$x_{-i}\in X_{-i},$\textit{\ is upper
semicontinuous;}

\textit{d) for each }$(x_{-i},y_{-i})\in X_{-i}\times X_{-i},$\textit{\ }$%
\{(x_{i},y_{i})\in X_{i}\times X_{i}:$\textit{\ }$f_{i}(x,y,u_{i})\nsubseteq 
$\textit{int}$C_{i}(x_{-i})$\textit{\ for each }$u_{i}\in A_{i}(x_{-i})\}$%
\textit{\ is nonempty}$;$

\textit{e) for each }$(x_{-i},y_{-i},u_{i})\in X_{-i}\times X_{-i}\times
X_{i},$\textit{\ }$f_{i}(x_{-i},\cdot ,y_{-i},\cdot ,u_{i}):X_{i}\times
X_{i}\rightarrow 2^{Y_{i}}$\textit{\ is int}$C_{i}(x_{-i})-$\textit{\
quasi-convex.}

\textit{Then, there exists a pair }$(x^{\ast },y^{\ast })\in X\times X$%
\textit{\ such that for each }$i\in N,$\textit{\ }$x_{i}^{\ast }\in 
\overline{A}_{i}(x_{-i}^{\ast }),$\textit{\ }$y_{i}^{\ast }\in \overline{F}%
_{i}(y_{-i}^{\ast })$\textit{\ and } $f_{i}(x^{\ast },y^{\ast
},u_{i})\nsubseteq $int$C_{i}(x_{-i}^{\ast })=\emptyset $\textit{\ for each }%
$u_{i}\in A_{i}(x_{-i}^{\ast }),$\textit{\ that is, }$(x^{\ast },y^{\ast })$%
\textit{\ is a solution for SGVQEP (IV).\medskip }

\section{Concluding remarks}

The first objective of this study has been to investigate the existence of
equilibrium for generalized games in choice form and to derive new, simpler
results on the existence of equilibrium for generalized abstract economies.
Then, our work has focused on applications concerning the solvability of
several types of systems of vector quasi-equilibrium problems. Our approach
has led to new and simple hypotheses which ensure the existence of solutions
for the classes of equilibrium problems we had considered.

\section{Appendix}

\textit{Proof of Theorem 1.} For each $i\in N,$ let us define the
correspondence $T_{i}:X_{-i}\times X_{-i}\rightarrow 2^{X_{i}\times X_{i}},$
by

$T_{i}(x_{-i},y_{-i})=\left\{ 
\begin{array}{c}
\text{co}(\bigcup\nolimits_{\{(z_{-i},t_{-i}):\mathcal{C}_{i}(z_{-i},t_{-i})%
\neq \emptyset \}}\mathcal{C}_{i}(z_{-i},t_{-i}))\text{ if }%
(x_{-i},y_{-i})\notin W_{i}; \\ 
\mathcal{C}_{i}(x_{-i},y_{-i})\text{ if }(x_{-i},y_{-i})\in W_{i}.%
\end{array}%
\right. $

The correspondence $T_{i}$ has nonempty and convex values$.$

If $(x_{i},y_{i})\in X_{i}\times X_{i},$ then, $T_{i}^{-1}(x_{i},y_{i})=%
\{(x_{-i},y_{-i})\in X_{-i}\times X_{-i}:(x_{i},y_{i})\in
T_{i}(x_{-i},y_{-i})\}=^{C}W_{i}\cup \mathcal{C}(x_{i},y_{i})$ is an open
set.

We apply the Yannelis and Prabhakar's Lemma and we obtain that $T_{i}$%
\textit{\ }has a continuous selection $f_{i}:X_{-i}\times X_{-i}\rightarrow
X_{i}\times X_{i}.$

Let $f:X\times X\rightarrow X\times X$ be defined by $f(x,y)=\prod%
\nolimits_{i\in N}f_{i}(x_{-i},y_{-i})$ for each $(x,y)\in X\times X.$ The
function $f$ is continuous, and, according to the Tychonoff fixed point
Theorem \cite{is}, there exists $(x^{\ast },y^{\ast })\in X\times X$ such
that $f(x^{\ast },y^{\ast })=(x^{\ast },y^{\ast }).$ Hence, $(x^{\ast
},y^{\ast })\in \prod\nolimits_{i\in N}T_{i}(x_{-i}^{\ast },y_{-i}^{\ast })$
and obviously, $(x_{i}^{\ast },y_{i}^{\ast })\in T_{i}(x_{-i}^{\ast
},y_{-i}^{\ast })$ for each $i\in N.$ Suppose that $((x_{-i}^{\ast
},A_{i}(x_{-i}^{\ast })),(y_{-i}^{\ast },F_{i}(y_{-i}^{\ast }))\cap \mathcal{%
C}_{i}\neq \emptyset ,$ for some $i\in N.$ Then, $\mathcal{C}%
_{i}(x_{-i}^{\ast },y_{-i}^{\ast })\neq \emptyset $ and $(x_{i}^{\ast
},y_{i}^{\ast })\in \mathcal{C}_{i}(x_{-i}^{\ast },y_{-i}^{\ast }),$ which
implies $(x^{\ast },y^{\ast })\in \mathcal{C}_{i}.\medskip \medskip $

\textit{Proof of Theorem 3. }For each $i\in N,$ let us define the set $%
\mathcal{C}_{i}=\{(x,y)\in X\times X:$ $A_{i}(x_{-i})\cap
P_{i}(x,y)=\emptyset \}\cap ($Gr$A_{i}\times $Gr$F_{i}).$

Then, $\mathcal{C}_{i}(x_{-i},y_{-i})=\{(x_{i},y_{i})\in X_{i}\times X_{i}:$ 
$A_{i}(x_{-i})\cap P_{i}(x,y)=\emptyset \}\cap (A_{i}(x_{-i})\times
F_{i}(y_{-i}))$ for each $(x_{-i},y_{-i})\in X_{-i}\times X_{-i}$ and

$\mathcal{C}_{i}(x_{i},y_{i})=(A_{i}^{-1}(x_{i})\times
F_{i}^{-1}(y_{i}))\cap \{(x_{-i},y_{-i})\in X_{-i}\times
X_{-i}:A_{i}(x_{-i})\cap P_{i}(x,y)=\emptyset \}\ $for each $%
(x_{i},y_{i})\in X_{i}\times X_{i}.$

Assumption c) implies that $\mathcal{C}_{i}$ is nonempty$.$ The set $%
\mathcal{C}_{i}(x_{i},y_{i})$ is open for each $(x_{i},y_{i})\in X_{i}\times
X_{i}$ since Assumptions b) and e) hold.

According to Assumptions b) and d), $(A_{i}(x_{-i})\times F_{i}(y_{-i}))\cap 
\mathcal{C}_{i}(x_{-i},y_{-i})\ $ is nonempty and convex for each $%
(x_{-i},y_{-i})\in X_{-i}\times X_{-i}.$

All hypotheses of Theorem 2 are fulfilled, and then, there exists $(x^{\ast
},y^{\ast })$ equilibrium in choice for the generalized game $%
(X_{i},A_{i},F_{i},\mathcal{C}_{i})_{i\in N}$. Obviously, $(x^{\ast
},y^{\ast })$ is equilibrium for the abstract economy $%
(X_{i},A_{i},F_{i},P_{i})_{i\in N}.\medskip $

\textit{Proof of Theorem 5}. For each $i\in N,$ let us define the set $%
\mathcal{C}_{i}=\{(x,y)\in X\times X:$ $A_{i}(x_{-i})\cap
P_{i}(x,y)=\emptyset \}\cap ($Gr$\overline{A_{i}}\times $Gr$\overline{F_{i}}%
).$

We claim that $\mathcal{C}_{i}$ is a closed set. Indeed, let us consider a
sequence $(x^{n},y^{n})_{n}$ in $\mathcal{C}_{i}$ and let us assume that $%
\lim_{n\rightarrow \infty }(x^{n},y^{n})=(x^{0},y^{0}).$ We will prove that $%
(x^{0},y^{0})\in \mathcal{C}_{i}.$

Since Gr$\overline{A_{i}}$ and Gr$\overline{F_{i}}$ are closed, it follows
that $(x^{0},y^{0})\in $Gr$\overline{A_{i}}\times $Gr$\overline{F_{i}}.$ It
remain to prove that $($ $x^{0},y^{0})\in \{(x,y)\in X\times X:$ $%
A_{i}(x_{-i})\cap P_{i}(x,y)=\emptyset \}.$

Suppose, by contrary, that the last assertion is false. Then, there exists $%
z_{i}^{0}\in X_{i}$ such that $z_{i}^{0}\in A_{i}(x_{-i}^{0})\cap
P_{i}(x^{0},y^{0}).$ The lower semicontinuity of $A_{i}$ implies the
existence of the sequence $(z_{i}^{n})_{n}$ in $X_{i}$ such that $%
\lim_{n\rightarrow \infty }z_{i}^{n}=z_{i}^{0}$ and $z_{i}^{n}\in
A_{i}(x_{-i}^{n})$ for each $n\in \mathbb{N}.$ We have also $z_{i}^{0}\in
P_{i}(x^{0},y^{0}),$ that is, $(x^{0},y^{0},z_{i}^{0})\in $Gr$P_{i}.$ Using
the fact that Gr$P_{i}$ is an open set, we can find open neighbourhoods $%
V_{(x^{0},y^{0})}$ and $V_{z_{i}^{0}}$ for $(x^{0},y^{0}),$ respectively $%
z_{i}^{0}$ such that $V_{(x^{0},y^{0})}\times V_{z_{i}^{0}}$ $\subseteq $Gr$%
P_{i}.$ It follows that $(x^{n},y^{n},z_{i}^{n})\in $Gr$P_{i}$ for all but
finitely many values of $n.$ Thus, $z_{i}^{n}\in P_{i}(x^{n},y^{n})$ for all
but finitely many values of $n.$ This contradicts the fact that the sequence 
$(x^{n},y^{n})_{n}$ $\ $\ is in $\mathcal{C}_{i}.$ Consequently, $($ $%
x^{0},y^{0})\in \{(x,y)\in X\times X:$ $A_{i}(x_{-i})\cap
P_{i}(x,y)=\emptyset \}$ and the set $\mathcal{C}_{i}$ is closed.

Let us define the correspondence $T_{i}:X_{-i}\times X_{-i}\rightarrow
2^{X_{i}\times X_{i}}$ by

$T_{i}(x_{-i},y_{-i})=\mathcal{C}_{i}(x_{-i},y_{-i})$ if $(x_{-i},y_{-i})\in
X_{-i}\times X_{-i}.$

We notice that $\mathcal{C}_{i}(x_{-i},y_{-i})=\{(x_{i},y_{i})\in
X_{i}\times X_{i}:$ $A_{i}(x_{-i})\cap P_{i}(x,y)=\emptyset \}\cap (%
\overline{A_{i}}(x_{-i})\times \overline{F_{i}}(y_{-i}))$ for each $%
(x_{-i},y_{-i})\in X_{-i}\times X_{-i}.$

Assumptions c) and e) imply that the correspondence $T_{i}$ has nonempty and
convex values$.$ $T_{i}$ is also closed valued.

Gr$T_{i}=\{(x,y)\in X\times X:$ $(x_{i},y_{i})\in \mathcal{C}%
_{i}(x_{-i},y_{-i})\}=\{(x,y)\in X\times X:$ $(x_{i},y_{i})\in \overline{%
A_{i}}(x_{-i})\times \overline{F_{i}}(y_{-i}),A_{i}(x_{-i})\cap
P_{i}(x,y)=\emptyset $ $\}=\mathcal{C}_{i}.$

Since Gr$T_{i}$ is closed and $X_{i}\times X_{i}$ is compact, by Theorem
7.1.16 from Klein and Thompson \cite{kl}, \ the correspondence $T_{i}$ is
upper semicontinuous.

Let $T:X\times X\rightarrow 2^{X\times X}$ be defined by $%
T(x,y)=\prod\nolimits_{i\in N}T_{i}(x_{-i},y_{-i})$ for each $(x,y)\in
X\times X.$ The correspondence $T$ is upper semicontinuous, and, according
to the Ky Fan fixed point Theorem \cite{fan}, there exists $(x^{\ast
},y^{\ast })\in X\times X$ such that $(x^{\ast },y^{\ast })\in T(x^{\ast
},y^{\ast }).$ Obviously, $(x_{i}^{\ast },y_{i}^{\ast })\in
T_{i}(x_{-i}^{\ast },y_{-i}^{\ast })$ for each $i\in N.$ Then, for each $%
i\in N,$ $\mathcal{C}_{i}(x_{-i}^{\ast },y_{-i}^{\ast })\neq \emptyset $ and 
$(x_{i}^{\ast },y_{i}^{\ast })\in \mathcal{C}_{i}(x_{-i}^{\ast
},y_{-i}^{\ast }),$ which implies $(x^{\ast },y^{\ast })\in \mathcal{C}%
_{i}.\medskip $ Therefore, for each $i\in N,$ $\ x_{i}^{\ast }\in \overline{%
A_{i}}(x_{-i}^{\ast }),$ $y_{i}^{\ast }\in \overline{F_{i}}(y_{-i}^{\ast })$
and $A_{i}(x_{-i}^{\ast })\cap P_{i}(x^{\ast },y^{\ast })=\emptyset
.\medskip $

\textit{Proof of Theorem 6.} For each $i\in N,$ let us define $P_{i}:X\times
X\rightarrow 2^{X_{i}},$ by

$P_{i}(x,y)=\{u_{i}\in X_{i}:$ $f_{i}(x,y,u_{i})\nsubseteq $int$%
C_{i}(x_{-i})\}$ for each $(x,y)\in X\times X.$

Assumption d) implies that $\{(x_{i},y_{i})\in A_{i}(x_{-i})\times
F_{i}(y_{-i}):\overline{A}_{i}(x_{-i})\cap P_{i}(x,y)=\emptyset \}\mathit{\ }
$ is nonempty for each $(x_{-i},y_{-i})\in X_{-i}\times X_{-i}.$

Let us consider $(x_{-i}^{0},y_{-i}^{0})\in X_{-i}\times X_{-i}.$ Now, we
shall prove that $D_{i}(x_{-i}^{0},y_{-i}^{0})$ is convex, where

$D_{i}(x_{-i}^{0},y_{-i}^{0})=\{(x_{i},y_{i})\in X_{i}\times X_{i}:\overline{%
A}_{i}(x_{-i}^{0})\cap P_{i}(x_{-i}^{0},x_{i},y_{-i}^{0},y_{i})=\emptyset \}%
\mathit{=}$

$=\{(x_{i},y_{i})\in X_{i}\times
X_{i}:f_{i}(x_{-i}^{0},x_{i},y_{-i}^{0},y_{i},u_{i})\subseteq $int$%
C_{i}(x_{-i}^{0})$ for each $u_{i}\in \overline{A}_{i}(x_{-i}^{0})\}.$

Let $(x_{i}^{1},y_{i}^{1}),(x_{i}^{2},y_{i}^{2})\in
D_{i}(x_{-i}^{0},y_{-i}^{0}).$ This means that:

$f_{i}(x_{-i}^{0},x_{i}^{1},y_{-i}^{0},y_{i}^{1},u_{i})\subseteq \medskip $%
int$C_{i}(x_{-i}^{0})$ for each $u_{i}\in \overline{A}_{i}(x_{-i}^{0})$ and

$f_{i}(x_{-i}^{0},x_{i}^{2},y_{-i}^{0},y_{i}^{2},u_{i})\subseteq $int$%
C_{i}(x_{-i}^{0})$ for each $u_{i}\in \overline{A}_{i}(x_{-i}^{0}).$

Let $\lambda \in \lbrack 0,1]$ and $(x_{i}(\lambda ),y_{i}(\lambda
))=\lambda (x_{i}^{1},y_{i}^{1})+(1-\lambda )(x_{i}^{2},y_{i}^{2}).$ We
claim that $(x_{i}(\lambda ),y_{i}(\lambda ))\in
D_{i}(x_{-i}^{0},y_{-i}^{0}).$

Since $f_{i}(x_{-i}^{0},\cdot ,y_{-i}^{0},\cdot ,u_{i}):X_{i}\times
X_{i}\rightarrow 2^{Y_{i}}$ is\textit{\ }$(-$int$C_{i}(x_{-i}^{0}))-$\textit{%
\ }quasi-convex-like for each $u_{i}\in X_{i},$ we have that $%
f_{i}(x_{-i}^{0},x_{i}(\lambda ),y_{-i}^{0},y_{i}(\lambda ),u_{i})\subseteq
f_{i}(x_{-i}^{0},x_{i}^{1},y_{-i}^{0},y_{i}^{1},u_{i})+$int$%
C_{i}(x_{-i}^{0}) $ or $f_{i}(x_{-i}^{0},x_{i}(\lambda
),y_{-i}^{0},y_{i}(\lambda ),u_{i})\subseteq
f_{i}(x_{-i}^{0},x_{i}^{2},y_{-i}^{0},y_{i}^{2},u_{i})+$int$%
C_{i}(x_{-i}^{0}).$ On the other hand, it is true that $%
f_{i}(x_{-i}^{0},x_{i}^{1},y_{-i}^{0},y_{i}^{1},u_{i})\subseteq $int$%
C_{i}(x_{-i}^{0})$ for each $u_{i}\in \overline{A}_{i}(x_{-i}^{0})$ and $%
f_{i}(x_{-i}^{0},x_{i}^{2},y_{-i}^{0},y_{i}^{2},u_{i})\subseteq $int$%
C_{i}(x_{-i}^{0})$ for each $u_{i}\in \overline{A}_{i}(x_{-i}^{0}).$ We
obtain that $f_{i}(x_{-i}^{0},x_{i}(\lambda ),y_{-i}^{0},y_{i}(\lambda
),u_{i})\subseteq $int$C_{i}(x_{-i}^{0})$ for each $u_{i}\in \overline{A}%
_{i}(x_{-i}^{0})$, that is, $(x_{i}(\lambda ),y_{i}(\lambda ))\in
D_{i}(x_{-i}^{0},y_{-i}^{0})$. Therefore, $D_{i}(x_{-i}^{0},y_{-i}^{0})$ is
convex.

Let $E_{i}(x_{i},y_{i})=\{(x_{-i},y_{-i})\in X_{-i}\times X_{-i}:$ $%
f_{i}(x,y,u_{i})\subseteq $int$C_{i}(x_{-i})$\textit{\ }for each\textit{\ }$%
u_{i}\in \overline{A}_{i}(x_{-i})\}$ be defined for each $(x_{i},y_{i})\in
X_{i}\times X_{i}.$ \newline
We claim that, for each $(x_{i},y_{i})\in X_{i}\times X_{i},$ $%
^{C}E_{i}(x_{i},y_{i})$ is closed, where $^{C}E_{i}(x_{i},y_{i})=$

$=\{(x_{-i},y_{-i})\in X_{-i}\times X_{-i}:$ $\exists u_{i}\in \overline{A}%
_{i}(x_{-i})$ such that $f_{i}(x,y,u_{i})\nsubseteq $int$C_{i}(x_{-i})\}$.

Indeed, let $(x_{i},y_{i})$ be fixed and let $(x_{-i}^{0},y_{-i}^{0})\in $cl$%
^{C}E_{i}(x_{i},y_{i}).$ Then, there exists $(x_{-i}^{n},y_{-i}^{n})_{n}$ a
sequence in $^{C}E_{i}(x_{i},y_{i})$ such that $\lim_{n\rightarrow \infty
}(x_{-i}^{n},y_{-i}^{n})=(x_{-i}^{0},y_{-i}^{0}).$ Then, for each $n\in N,$ $%
\exists u_{i}^{n}\in \overline{A}_{i}(x_{-i}^{n})$ such that $%
f_{i}(x_{-i}^{n},x_{i},y_{-i}^{n},y_{i}^{n},u_{i}^{n})\nsubseteq $int$%
C_{i}(x_{-i}^{n}),$ that is, $%
f_{i}(x_{-i}^{n},x_{i},y_{-i}^{n},y_{i},u_{i}^{n})\cap W_{i}(x_{-i}^{n})\neq
\emptyset ,$ where $W_{i}(x_{-i}^{n})=Y_{i}\backslash (-$int$%
C_{i}(x_{-i}^{n})).$ The set $X_{i}$ is compact and therefore, we can
assume, without generality, that the sequence $(u_{i}^{n})_{n}$ is
convergent and let $\lim_{n\rightarrow \infty }u_{i}^{n}=u_{i}^{0}\in
X_{-i}. $ The closedness of $\overline{A}_{i}$ implies that $u_{i}^{0}\in 
\overline{A}_{i}(x_{-i}^{0}).$

We assert that there exists a sequence $(z_{i}^{n})_{n}$ in $X_{i}$ such
that $z_{i}^{n}\in f_{i}(x_{-i}^{n},x_{i},y_{-i}^{n},y_{i},u_{i}^{n})\cap
W_{i}(x_{-i}^{n})$ for each $n\in \mathbb{N}.$ It follows that $z_{i}^{n}\in
W_{i}(x_{-i}^{n})$ for each $n\in \mathbb{N}.$ Since $X_{i}$ is compact, we
can suppose that $\lim_{n\rightarrow \infty }z_{i}^{n}=z_{i}^{0}.$The
closedness of $W_{i}$ implies that $z_{i}^{0}\in W_{i}(x_{-i}^{0}).$ We
invoke here the closedness of $f_{i}(\cdot ,x_{i},\cdot ,y_{i},\cdot ):$%
\textit{\ }$X_{-i}\times X_{-i}\times X_{i}\rightarrow 2^{Y_{i}}$ and we
conclude that $z_{i}^{0}\in
f_{i}(x_{-i}^{0},x_{i},y_{-i}^{0},y_{i},u_{i}^{0}).$ Therefore, $%
f_{i}(x_{-i}^{0},x_{i},y_{-i}^{0},y_{i},u_{i}^{0})\cap W_{i}(x_{-i}^{0})\neq
\emptyset ,$ and, thus, $(x_{-i}^{0},y_{-i}^{0})\in ^{C}E_{i}(x_{i},y_{i}),$ 
$^{C}E_{i}(x_{i},y_{i})$ is closed and then, $E_{i}(x_{i},y_{i})=%
\{(x_{-i},y_{-i})\in X_{-i}\times X_{-i}:\overline{A}_{i}(x_{-i})\cap
P_{i}(x,y)=\emptyset \}$ is an open set.

All assumptions of Theorem 4 are fulfilled and there exists a pair $(x^{\ast
},y^{\ast })\in $\ $X\times X$\ such that $\ x_{i}^{\ast }\in
A_{i}(x_{-i}^{\ast }),$\ $y_{i}^{\ast }\in F_{i}(y_{-i}^{\ast })$\ and $%
A_{i}(x_{-i}^{\ast })\cap P_{i}(x^{\ast },y^{\ast })=\emptyset $\ for each $%
i\in N.\medskip $ Then, for each $i\in N,$ $x_{i}^{\ast }\in
A_{i}(x_{-i}^{\ast }),$ $y_{i}^{\ast }\in F_{i}(y_{-i}^{\ast })$ and $%
f_{i}(x^{\ast },y^{\ast },u_{i})\subseteq C_{i}(x_{-i}^{\ast })$ for each $%
u_{i}\in A_{i}(x_{-i}^{\ast }),$ that is, $(x^{\ast },y^{\ast })$ is a
solution for SGVQEP (I).\medskip 

\textit{Proof of Theorem 7. \ }For each $i\in N,$ let us define $%
P_{i}:X\times X\rightarrow 2^{X_{i}},$ by

$P_{i}(x,y)=\{u_{i}\in X_{i}:$ $f_{i}(x,y,u_{i})\nsubseteq C_{i}(x_{-i})\}$
for each $(x,y)\in X\times X.$

Assumption d) implies that $\{(x_{i},y_{i})\in X_{i}\times
X_{i}:A_{i}(x_{-i})\cap P_{i}(x,y)=\emptyset \}\mathit{\ }$ is nonempty for
each $(x_{-i},y_{-i})\in X_{-i}\times X_{-i}.$

Let us consider $(x_{-i}^{0},y_{-i}^{0})\in X_{-i}\times X_{-i}$ and $%
D_{i}(x_{-i}^{0},y_{-i}^{0})=\{(x_{i},y_{i})\in X_{i}\times
X_{i}:A_{i}(x_{-i}^{0})\cap
P_{i}(x_{-i}^{0},x_{i},y_{-i}^{0},y_{i})=\emptyset \}\mathit{.}$

We can prove that $D_{i}(x_{-i}^{0},y_{-i}^{0})$ is convex, as in the proof
of Theorem 6.

Following the same line as in the proof of Theorem 3.2.1. from \cite{lin2007}%
, we can show that $P_{i}$ has an open graph.

All assumptions of Theorem 5 are verified. According to this result, there
exists $(x^{\ast },y^{\ast })$ such that for each $i\in N,$ $\ x_{i}^{\ast
}\in \overline{A_{i}}(x_{-i}^{\ast }),$ $y_{i}^{\ast }\in \overline{F_{i}}%
(y_{-i}^{\ast })$ and $A_{i}(x_{-i}^{\ast })\cap P_{i}(x^{\ast },y^{\ast
})=\emptyset .$ Then, for each $i\in N,$ $x_{i}^{\ast }\in \overline{A}%
_{i}(x_{-i}^{\ast }),$ $y_{i}^{\ast }\in \overline{F}_{i}(y_{-i}^{\ast })$
and $f_{i}(x^{\ast },y^{\ast },u_{i})\subseteq C_{i}(x_{-i}^{\ast })$ for
each $u_{i}\in A_{i}(x_{-i}^{\ast }),$ that is, $(x^{\ast },y^{\ast })$ is a
solution for SGVQEP (I).\medskip 

\textit{Proof of Theorem 8. }For each $i\in N,$ let us define $P_{i}:X\times
X\rightarrow 2^{X_{i}},$ by

$P_{i}(x,y)=\{u_{i}\in X_{i}:$ $f_{i}(x,y,u_{i})\cap $int$%
C_{i}(x_{-i}))=\emptyset \}$ for each $(x,y)\in X\times X.$

Let $E_{i}(x_{i},y_{i})=\{(x_{-i},y_{-i})\in X_{-i}\times X_{-i}:$ $%
f_{i}(x,y,u_{i})\cap $int$C_{i}(x_{-i})\neq \emptyset $\textit{\ }for each%
\textit{\ }$u_{i}\in \overline{A}_{i}(x_{-i})\}$ be defined for each $%
(x_{i},y_{i})\in X_{i}\times X_{i}.$ We claim that, for each $%
(x_{i},y_{i})\in X_{i}\times X_{i}$, $^{C}E_{i}(x_{i},y_{i})$ is closed,
where $^{C}E_{i}(x_{i},y_{i})=\{(x_{-i},y_{-i})\in X_{-i}\times X_{-i}:$ $%
\exists u_{i}\in \overline{A}_{i}(x_{-i})$ such that $f_{i}(x,y,u_{i})\cap $%
int$C_{i}(x_{-i})=\emptyset \}$. Indeed, let $(x_{i},y_{i})$ be fixed and
let $(x_{-i}^{0},y_{-i}^{0})\in $cl$^{C}E_{i}(x_{i},y_{i}).$ Then, there
exists $(x_{-i}^{n},y_{-i}^{n})_{n}$ a sequence in $^{C}E_{i}(x_{i},y_{i})$
such that $\lim_{n\rightarrow \infty
}(x_{-i}^{n},y_{-i}^{n})=(x_{-i}^{0},y_{-i}^{0}).$Then, for each $n\in N,$ $%
\exists u_{i}^{n}\in \overline{A}_{i}(x_{-i}^{n})$ such that $%
f_{i}(x_{-i}^{n},x_{i},y_{-i}^{n},y_{i}^{n},u_{i}^{n})\cap $int$%
C_{i}(x_{-i}^{n})=\emptyset .$ The set $X_{i}$ is compact and therefore, we
can assume, without generality, that the sequence $(u_{i}^{n})_{n}$ is
convergent and let $\lim_{n\rightarrow \infty }u_{i}^{n}=u_{i}^{0}\in
X_{-i}. $ The closedness of $\overline{A}_{i}$ implies that $u_{i}^{0}\in 
\overline{A}_{i}(x_{-i}^{0}).$

Let $z_{i}^{0}\in f_{i}(x_{-i}^{0},x_{i},y_{-i}^{0},y_{i},u_{i}^{0}).$ We
invoke here the lower semicontinuity of $f_{i}(\cdot ,x_{i},\cdot
,y_{i},\cdot ):$\textit{\ }$X_{-i}\times X_{-i}\times X_{i}\rightarrow
2^{Y_{i}}$ and we assert that there exists a sequence $(z_{i}^{n})_{n}$ in $%
X_{i}$ such that $\lim_{n\rightarrow \infty }z_{i}^{n}=z_{i}^{0}$ and $%
z_{i}^{n}\in f_{i}(x_{-i}^{n},x_{i},y_{-i}^{n},y_{i},u_{i}^{n})\subseteq
W_{i}(x_{-i}^{n})$ for each $n\in \mathbb{N}.$ It follows that $z_{i}^{n}\in
W_{i}(x_{-i}^{n})$ for each $n\in \mathbb{N}.$ The closedness of $W_{i}$
implies that $z_{i}^{0}\in W_{i}(x_{-i}^{0}).$ Consequently, $%
f_{i}(x_{-i}^{0},x_{i},y_{-i}^{0},y_{i},u_{i}^{0})\subseteq
W_{i}(x_{-i}^{0}),$ and, thus, $(x_{-i}^{0},y_{-i}^{0})\in
^{C}E_{i}(x_{i},y_{i}),$ $^{C}E_{i}(x_{i},y_{i})$ is closed and then, $%
E_{i}(x_{i},y_{i})=\{(x_{-i},y_{-i})\in X_{-i}\times X_{-i}:\overline{A}%
_{i}(x_{-i})\cap P_{i}(x,y)=\emptyset \}$ is an open set.

Let us consider $(x_{-i}^{0},y_{-i}^{0})\in X_{-i}\times X_{-i}.$ Now, we
shall prove that $D(x_{-i}^{0},y_{-i}^{0})$ is convex, where

$D(x_{-i}^{0},y_{-i}^{0})=\{(x_{i},y_{i})\in X_{i}\times X_{i}:\overline{A}%
_{i}(x_{-i}^{0})\cap P_{i}(x_{-i}^{0},x_{i},y_{-i}^{0},y_{i})=\emptyset \}%
\mathit{=}$\newline
$=\{(x_{i},y_{i})\in X_{i}\times
X_{i}:f_{i}(x_{-i}^{0},x_{i},y_{-i}^{0},y_{i},u_{i})\cap $int$%
C_{i}(x_{-i}^{0})\neq \emptyset $ for each $u_{i}\in \overline{A}%
_{i}(x_{-i}^{0})\}.$\newline
Let $(x_{i}^{1},y_{i}^{1}),(x_{i}^{2},y_{i}^{2})\in
D(x_{-i}^{0},y_{-i}^{0}). $ This means $%
f_{i}(x_{-i}^{0},x_{i}^{1},y_{-i}^{0},y_{i}^{1},u_{i})\cap $int$%
C_{i}(x_{-i}^{0})\neq \emptyset $ for each $u_{i}\in \overline{A}%
_{i}(x_{-i}^{0})$ and $%
f_{i}(x_{-i}^{0},x_{i}^{2},y_{-i}^{0},y_{i}^{2},u_{i})\cap $int$%
C_{i}(x_{-i}^{0})\neq \emptyset $ for each $u_{i}\in \overline{A}%
_{i}(x_{-i}^{0}).$

For each $u_{i}\in \overline{A}_{i}(x_{-i}^{0}),$ let $z_{i}^{j}(u_{i})\in $ 
$f_{i}(x_{-i}^{0},x_{i}^{j},y_{-i}^{0},y_{i}^{j},u_{i})\cap $int$%
C_{i}(x_{-i}^{0})$ for each $j\in \{1,2\}.$ Since int$C_{i}(x_{-i}^{0})$ is
convex, it is true that co$\{z_{i}^{1}(u_{i}),z_{i}^{2}(u_{i})\}\subseteq $%
int$C_{i}(x_{-i}^{0})$ for each $u_{i}\in \overline{A}_{i}(x_{-i}^{0}).$ Let 
$(x_{i}(\lambda ),y_{i}(\lambda ))=\lambda (x_{i}^{1},y_{i}^{1})+(1-\lambda
)(x_{i}^{2},y_{i}^{2})$ be defined for each $\lambda \in \lbrack 0,1].$
Since for each $(x_{-i},y_{-i},u_{i})\in X_{-i}\times X_{-i}\times X_{i},$ $%
f_{i}(x_{-i},\cdot ,y_{-i},\cdot ,u_{i}):X_{i}\times X_{i}\rightarrow
2^{Y_{i}}$ is natural int$C_{i}(x_{-i})-$\ quasi-concave, it follows that $%
(x_{i}(\lambda ),y_{i}(\lambda ))\in D(x_{-i}^{0},y_{-i}^{0})$ for each $%
\lambda \in \lbrack 0,1].$

All assumptions of Theorem 4 are verified and there exists $(x^{\ast
},y^{\ast })$ such that for each $i\in N,$ $\ x_{i}^{\ast }\in
A_{i}(x_{-i}^{\ast }),$ $y_{i}^{\ast }\in F_{i}(y_{-i}^{\ast })$ and $%
A_{i}(x_{-i}^{\ast })\cap P_{i}(x^{\ast },y^{\ast })=\emptyset .$ Then, for
each $i\in N,$ $x_{i}^{\ast }\in \overline{A}_{i}(x_{-i}^{\ast }),$ $%
y_{i}^{\ast }\in \overline{F}_{i}(x_{-i}^{\ast })$ and $f_{i}(x^{\ast
},y^{\ast },u_{i})\cap C_{i}(x_{-i})\neq \emptyset $ for each $u_{i}\in
A_{i}(x_{-i}^{\ast }),$ that is, $(x^{\ast },y^{\ast })$ is a solution for
SGVQEP (II).\medskip

\textit{Proof of Theorem 9. }For each $i\in N,$ let us define $P_{i}:X\times
X\rightarrow 2^{X_{i}},$ by

$P_{i}(x,y)=\{u_{i}\in X_{i}:$ $f_{i}(x,y,u_{i})\cap C_{i}(x_{-i})=\emptyset
\}$ for each $(x,y)\in X\times X.$

Assumption d) implies that $\{(x_{i},y_{i})\in X_{i}\times
X_{i}:A_{i}(x_{-i})\cap P_{i}(x,y)=\emptyset \}\mathit{\ }$ is nonempty for
each $(x_{-i},y_{-i})\in X_{-i}\times X_{-i}.$

Let us consider $(x_{-i}^{0},y_{-i}^{0})\in X_{-i}\times X_{-i}.$ We can
prove that $D_{i}(x_{-i}^{0},y_{-i}^{0})$ is convex, by using a similar
argument as in the proof of Theorem 8, where $D_{i}(x_{-i}^{0},y_{-i}^{0})=%
\{(x_{i},y_{i})\in X_{i}\times X_{i}:A_{i}(x_{-i}^{0})\cap
P_{i}(x_{-i}^{0},x_{i},y_{-i}^{0},y_{i})=\emptyset \}.$

Now, we shall prove that $P_{i}$ has an open graph. In order to do this, let
us consider $(x^{0},y^{0},u_{i}^{0})\in $cl$^{C}$[Gr$P_{i}$]. We claim that $%
(x^{0},y^{0},u_{i}^{0})\in ^{C}$[Gr$P_{i}$] and, thus, $^{C}$[Gr$P_{i}$] is
closed, which implies that Gr$P_{i}$ is open. Since $(x^{0},y^{0},u_{i}^{0})%
\in $cl$^{C}$[Gr$P_{i}$], there exists a sequence $%
(x^{n},y^{n},u_{i}^{n})_{n}$ in cl$^{C}$[Gr$P_{i}$] such that $%
\lim_{n\rightarrow \infty }(x^{n},y^{n},u_{i}^{n})=(x^{0},y^{0},u_{i}^{0}).$ 
$\ $Therefore, $f_{i}(x^{n},y^{n},u_{i}^{n})\cap C_{i}(x_{-i}^{n})\neq
\emptyset $ for each $n\in \mathbb{N}.$ Let us consider $z_{i}^{n}\in
f_{i}(x^{n},y^{n},u_{i}^{n})\cap C_{i}(x_{-i}^{n}).$ We can assume that $%
\lim_{n\rightarrow \infty }z_{i}^{n}=z_{i}^{0}$. We invoke here the upper
semicontinuity of $f_{i}$ and $C_{i}$ and we conclude that $z_{i}^{0}\in $.$%
f_{i}(x^{0},y^{0},u_{i}^{0})\cap C_{i}(x_{-i}^{0})$ and, thus, $%
(x^{0},y^{0},u_{i}^{0})\in ^{C}$[Gr$P_{i}$].

All assumptions of Theorem 5 are verified. According to this result, there
exists $(x^{\ast },y^{\ast })$ such that for each $i\in N,$ $\ x_{i}^{\ast
}\in \overline{A}_{i}(x_{-i}^{\ast }),$ $y_{i}^{\ast }\in \overline{F}%
_{i}(y_{-i}^{\ast })$ and $A_{i}(x_{-i}^{\ast })\cap P_{i}(x^{\ast },y^{\ast
})=\emptyset .$ Then, for each $i\in N,$ $x_{i}^{\ast }\in \overline{A}%
_{i}(x_{-i}^{\ast }),$ $y_{i}^{\ast }\in \overline{F}_{i}(y_{-i}^{\ast })$
and $f_{i}(x^{\ast },y^{\ast },u_{i})\cap C_{i}(x_{-i})\neq \emptyset $ for
each $u_{i}\in A_{i}(x_{-i}^{\ast }),$ that is, $(x^{\ast },y^{\ast })$ is a
solution for SGVQEP (II).\medskip 

\textit{Proof of Theorem 10. }For each $i\in N,$ let us define $%
P_{i}:X\times X\rightarrow 2^{X_{i}},$ by

$P_{i}(x,y)=\{u_{i}\in X_{i}:$ $f_{i}(x,y,u_{i})\cap (-C_{i}(x_{-i}))\neq
\emptyset \}$ for each $(x,y)\in X\times X.$

Let $E_{i}(x_{i},y_{i})=\{(x_{-i},y_{-i})\in X_{-i}\times X_{-i}:$ $%
f_{i}(x,y,u_{i})\cap (-C_{i}(x_{-i}))=\emptyset $\textit{\ }for each\textit{%
\ }$u_{i}\in \overline{A}_{i}(x_{-i})\}$ be defined for each $%
(x_{i},y_{i})\in X_{i}\times X_{i}.$ We claim that, for each $%
(x_{i},y_{i})\in X_{i}\times X_{i},$ $^{C}E_{i}(x_{i},y_{i})$ is closed,
where $^{C}E_{i}(x_{i},y_{i})=\{(x_{-i},y_{-i})\in X_{-i}\times X_{-i}:$ $%
\exists u_{i}\in \overline{A}_{i}(x_{-i})$ such that $f_{i}(x,y,u_{i})\cap
(-C_{i}(x_{-i}))\neq \emptyset \}$.

Indeed, let $(x_{i},y_{i})$ be fixed and let $(x_{-i}^{0},y_{-i}^{0})\in $cl$%
^{C}E_{i}(x_{i},y_{i}).$ Then, there exists $(x_{-i}^{n},y_{-i}^{n})_{n}$ a
sequence in $^{C}E_{i}(x_{i},y_{i})$ such that $\lim_{n\rightarrow \infty
}(x_{-i}^{n},y_{-i}^{n})=(x_{-i}^{0},y_{-i}^{0}).$ Then, for each $n\in N,$ $%
\exists u_{i}^{n}\in \overline{A}_{i}(x_{-i}^{n})$ such that $%
f_{i}(x_{-i}^{n},x_{i},y_{-i}^{n},y_{i}^{n},u_{i}^{n})\cap
(-C_{i}(x_{-i}^{n}))\neq \emptyset .$ The set $X_{i}$ is compact and
therefore, we can assume, without generality, that the sequence $%
(u_{i}^{n})_{n}$ is convergent and let $\lim_{n\rightarrow \infty
}u_{i}^{n}=u_{i}^{0}\in X_{-i}.$ The closedness of $\overline{A}_{i}$
implies that $u_{i}^{0}\in \overline{A}_{i}(x_{-i}^{0}).$

Let $z_{i}^{n}\in f_{i}(x_{-i}^{n},x_{i},y_{-i}^{n},y_{i},u_{i}^{n})\cap
(-C_{i}(x_{-i}^{n})).$ We invoke here the upper semicontinuity of $%
f_{i}(\cdot ,x_{i},\cdot ,y_{i},\cdot ):$\textit{\ }$X_{-i}\times
X_{-i}\times X_{i}\rightarrow 2^{Y_{i}}$ and we assert that $%
\lim_{n\rightarrow \infty }z_{i}^{n}=z_{i}^{0}.$ It follows that $%
z_{i}^{n}\in -C_{i}(x_{-i}^{n})$ for each $n\in \mathbb{N}.$ The closedness
of $C_{i}$ implies that $z_{i}^{0}\in $ $-C_{i}(x_{-i}^{0})$ and the the
closedness of $f_{i}(x_{-i},\cdot ,y_{-i},\cdot ,\cdot )$ implies that $%
z_{i}^{0}\in $ $f_{i}(x_{-i}^{0},x_{i},y_{-i}^{0},y_{i},u_{i}^{0}).$
Consequently, there exists $u_{i}^{0}\in \overline{A}_{i}(x_{-i}^{0})$ such
that $f_{i}(x_{-i}^{0},x_{i},y_{-i}^{0},y_{i},u_{i}^{0})\cap
(-C_{i}(x_{-i}^{0}))\neq \emptyset ,$ and, thus, $(x_{-i}^{0},y_{-i}^{0})\in
^{C}E_{i}(x_{i},y_{i}),$ $^{C}E_{i}(x_{i},y_{i})$ is closed and then, $%
E_{i}(x_{i},y_{i})=\{(x_{-i},y_{-i})\in X_{-i}\times X_{-i}:\overline{A}%
_{i}(x_{-i})\cap P_{i}(x,y)=\emptyset \}$ is an open set.

Let us consider $(x_{-i}^{0},y_{-i}^{0})\in X_{-i}\times X_{-i}.$ Now, we
shall prove that $D(x_{-i}^{0},y_{-i}^{0})$ is convex,where

$D(x_{-i}^{0},y_{-i}^{0})=\{(x_{i},y_{i})\in X_{i}\times X_{i}:\overline{A}%
_{i}(x_{-i}^{0})\cap P_{i}(x_{-i}^{0},x_{i},y_{-i}^{0},y_{i})=\emptyset \}%
\mathit{=}$\newline
$=\{(x_{i},y_{i})\in X_{i}\times
X_{i}:f_{i}(x_{-i}^{0},x_{i},y_{-i}^{0},y_{i},u_{i})\cap ($-$%
C_{i}(x_{-i}^{0}))\neq \emptyset $ for each $u_{i}\in \overline{A}%
_{i}(x_{-i}^{0})\}.$

Let $(x_{i}^{1},y_{i}^{1}),(x_{i}^{2},y_{i}^{2})\in
D(x_{-i}^{0},y_{-i}^{0}). $ This means:

$f_{i}(x_{-i}^{0},x_{i}^{1},y_{-i}^{0},y_{i}^{1},u_{i})\cap ($-$%
C_{i}(x_{-i}^{0}))=\emptyset $ for each $u_{i}\in \overline{A}%
_{i}(x_{-i}^{0})$ and

$f_{i}(x_{-i}^{0},x_{i}^{2},y_{-i}^{0},y_{i}^{2},u_{i})\cap ($-$%
C_{i}(x_{-i}^{0}))=\emptyset $ for each $u_{i}\in \overline{A}%
_{i}(x_{-i}^{0}).$

Let $(x_{i}(\lambda ),y_{i}(\lambda ))=\lambda
(x_{i}^{1},y_{i}^{1})+(1-\lambda )(x_{i}^{2},y_{i}^{2})$ be defined for each 
$\lambda \in \lbrack 0,1].$ We claim that $(x_{i}(\lambda ),y_{i}(\lambda
))\in D(x_{-i}^{0},y_{-i}^{0})$ for each $\lambda \in \lbrack 0,1].$

Indeed, let fix $\lambda \in \lbrack 0,1].$ Since $f_{i}(x_{-i}^{0},\cdot
,y_{-i}^{0},\cdot ,u_{i}):X_{i}\times X_{i}\rightarrow 2^{Y_{i}}$ is\textit{%
\ }$(-C_{i}(x_{-i}^{0}))-$\textit{\ }quasi-convex-like for each $u_{i}\in
X_{i},$ we have that:

$f_{i}(x_{-i}^{0},x_{i}(\lambda ),y_{-i}^{0},y_{i}(\lambda ),u_{i})\subseteq
f_{i}(x_{-i}^{0},x_{i}^{1},y_{-i}^{0},y_{i}^{1},u_{i})-(-C_{i}(x_{-i}^{0}))$
or

$f_{i}(x_{-i}^{0},x_{i}(\lambda ),y_{-i}^{0},y_{i}(\lambda ),u_{i})\subseteq
f_{i}(x_{-i}^{0},x_{i}^{2},y_{-i}^{0},y_{i}^{2},u_{i})-(-C_{i}(x_{-i}^{0})).$

On the other hand, it is true that $%
f_{i}(x_{-i}^{0},x_{i}^{1},y_{-i}^{0},y_{i}^{1},u_{i})\cap ($-$%
C_{i}(x_{-i}^{0}))=\emptyset $ for each $u_{i}\in \overline{A}%
_{i}(x_{-i}^{0})$ and $%
f_{i}(x_{-i}^{0},x_{i}^{2},y_{-i}^{0},y_{i}^{2},u_{i})\cap ($-$%
C_{i}(x_{-i}^{0}))=\emptyset $ for each $u_{i}\in \overline{A}%
_{i}(x_{-i}^{0}).$ We obtain that $f_{i}(x_{-i}^{0},x_{i}(\lambda
),y_{-i}^{0},y_{i}(\lambda ),u_{i})\cap ($-$C_{i}(x_{-i}^{0}))=\emptyset $
for each $u_{i}\in \overline{A}_{i}(x_{-i}^{0})$, that is, $(x_{i}(\lambda
),y_{i}(\lambda ))\in D(x_{-i}^{0},y_{-i}^{0})$. Consequently, $%
D(x_{-i}^{0},y_{-i}^{0})$ is convex.

All assumptions of Theorem 4 are fulfilled and there exists $(x^{\ast
},y^{\ast })$ an equilibrium for the associated generalized abstract economy 
$(X_{i},A_{i},F_{i},P_{i})_{i\in N}$. Then, for each $i\in N,$ $x_{i}^{\ast
}\in A_{i}(x_{-i}^{\ast }),$ $y_{i}^{\ast }\in F_{i}(y_{-i}^{\ast })$ and $%
f_{i}(x^{\ast },y^{\ast },u_{i})\cap (-$int$C_{i}(x_{-i}^{\ast }))=\emptyset 
$ for each $u_{i}\in A_{i}(x_{-i}^{\ast }),$ that is, $(x^{\ast },y^{\ast })$
is a solution for SGVQEP (III).\medskip 

\textit{Proof of Theorem 11. }For each $i\in N,$ let us define $%
P_{i}:X\times X\rightarrow 2^{X_{i}},$ by

$P_{i}(x,y)=\{u_{i}\in X_{i}:$ $f_{i}(x,y,u_{i})\cap (-$int$%
C_{i}(x_{-i}))\neq \emptyset \}$ for each $(x,y)\in X\times X.$

Assumption d) implies that $\{(x_{i},y_{i})\in X_{i}\times
X_{i}:A_{i}(x_{-i})\cap P_{i}(x,y)=\emptyset \}\mathit{\ }$ is nonempty for
each $(x_{-i},y_{-i})\in X_{-i}\times X_{-i}.$

Let us consider $(x_{-i}^{0},y_{-i}^{0})\in X_{-i}\times X_{-i}.$ Now, as in
the proof of Theorem 10, we can show that $D_{i}(x_{-i}^{0},y_{-i}^{0})$ is
convex$,$where

$D_{i}(x_{-i}^{0},y_{-i}^{0})=\{(x_{i},y_{i})\in X_{i}\times
X_{i}:A_{i}(x_{-i}^{0})\cap
P_{i}(x_{-i}^{0},x_{i},y_{-i}^{0},y_{i})=\emptyset \}.$

By following a similar argument as in the proof of Theorem 7, we obtain that
Gr$P_{i}$ is open.

All assumptions of Theorem 5 are verified. According to this result, there
exists $(x^{\ast },y^{\ast })$ an equilibrium for the associated generalized
abstract economy $(X_{i},A_{i},F_{i},P_{i})_{i\in N}.$ Then, for each $i\in
N,$ $x_{i}^{\ast }\in \overline{A}_{i}(x_{-i}^{\ast }),$ $y_{i}^{\ast }\in 
\overline{F}_{i}(x_{-i}^{\ast })$ and $f_{i}(x^{\ast },y^{\ast },u_{i})\cap
(-$int$C_{i}(x_{-i}^{\ast }))=\emptyset $ for each $u_{i}\in
A_{i}(x_{-i}^{\ast }),$ that is, $(x^{\ast },y^{\ast })$ is a solution for
SGVQEP (III).\medskip

\textit{Proof of Theorem 12.} For each $i\in N,$ let us define $%
P_{i}:X\times X\rightarrow 2^{X_{i}},$ by

$P_{i}(x,y)=\{u_{i}\in X_{i}:$ $f_{i}(x,y,u_{i})\subseteq C_{i}(x_{-i})\}$
for each $(x,y)\in X\times X.$

Let $E_{i}(x_{i},y_{i})=\{(x_{-i},y_{-i})\in X_{-i}\times X_{-i}:$ $%
f_{i}(x,y,u_{i})\nsubseteq C_{i}(x_{-i})$\textit{\ }for each\textit{\ }$%
u_{i}\in \overline{A}_{i}(x_{-i})\}$ be defined for each $(x_{i},y_{i})\in
X_{i}\times X_{i}.$ We claim that, for each $(x_{i},y_{i})\in X_{i}\times
X_{i}$, $^{C}E_{i}(x_{i},y_{i})$ is closed, where $^{C}E_{i}(x_{i},y_{i})=%
\{(x_{-i},y_{-i})\in X_{-i}\times X_{-i}:$ $\exists u_{i}\in \overline{A}%
_{i}(x_{-i})$ such that $f_{i}(x,y,u_{i})\subseteq C_{i}(x_{-i})\}$.

Indeed, let $(x_{i},y_{i})$ be fixed and let $(x_{-i}^{0},y_{-i}^{0})\in $cl$%
^{C}E_{i}(x_{i},y_{i}).$ Then, there exists $(x_{-i}^{n},y_{-i}^{n})_{n}$ a
sequence in $^{C}E_{i}(x_{i},y_{i})$ such that $\lim_{n\rightarrow \infty
}(x_{-i}^{n},y_{-i}^{n})=(x_{-i}^{0},y_{-i}^{0}).$Then, for each $n\in N,$ $%
\exists u_{i}^{n}\in \overline{A}_{i}(x_{-i}^{n})$ such that $%
f_{i}(x_{-i}^{n},x_{i},y_{-i}^{n},y_{i}^{n},u_{i}^{n})\subseteq
C_{i}(x_{-i}^{n}).$ The set $X_{i}$ is compact and therefore, we can assume,
without generality, that the sequence $(u_{i}^{n})_{n}$ is convergent and
let $\lim_{n\rightarrow \infty }u_{i}^{n}=u_{i}^{0}\in X_{-i}.$ The
closedness of $\overline{A}_{i}$ implies that $u_{i}^{0}\in \overline{A}%
_{i}(x_{-i}^{0}).$

Let $z_{i}^{0}\in f_{i}(x_{-i}^{0},x_{i},y_{-i}^{0},y_{i},u_{i}^{0}).$ We
invoke here the lower semicontinuity of $f_{i}(\cdot ,x_{i},\cdot
,y_{i},\cdot ):$\textit{\ }$X_{-i}\times X_{-i}\times X_{i}\rightarrow
2^{Y_{i}}$ and we assert that there exists a sequence $(z_{i}^{n})_{n}$ in $%
X_{i}$ such that $\lim_{n\rightarrow \infty }z_{i}^{n}=z_{i}^{0}$ and $%
z_{i}^{n}\in f_{i}(x_{-i}^{n},x_{i},y_{-i}^{n},y_{i},u_{i}^{n})\subseteq
C_{i}(x_{-i}^{n})$ for each $n\in \mathbb{N}.$ It follows that $z_{i}^{n}\in
C_{i}(x_{-i}^{n})$ for each $n\in \mathbb{N}.$ The closedness of $C_{i}$
implies that $z_{i}^{0}\in C_{i}(x_{-i}^{0}).$ Consequently, $%
f_{i}(x_{-i}^{0},x_{i},y_{-i}^{0},y_{i},u_{i}^{0})\subseteq
C_{i}(x_{-i}^{0}),$ and, thus, $(x_{-i}^{0},y_{-i}^{0})\in
^{C}E_{i}(x_{i},y_{i}),$ $^{C}E_{i}(x_{i},y_{i})$ is closed and then, $%
E_{i}(x_{i},y_{i})=\{(x_{-i},y_{-i})\in X_{-i}\times X_{-i}:\overline{A}%
_{i}(x_{-i})\cap P_{i}(x,y)=\emptyset \}$ is an open set.

Let us consider $(x_{-i}^{0},y_{-i}^{0})\in X_{-i}\times X_{-i}.$ Now, we
shall prove that $D(x_{-i}^{0},y_{-i}^{0})$ is convex, where

$D(x_{-i}^{0},y_{-i}^{0})=\{(x_{i},y_{i})\in X_{i}\times X_{i}:\overline{A}%
_{i}(x_{-i}^{0})\cap P_{i}(x_{-i}^{0},x_{i},y_{-i}^{0},y_{i})=\emptyset \}%
\mathit{=}$\newline
$=\{(x_{i},y_{i})\in X_{i}\times
X_{i}:f_{i}(x_{-i}^{0},x_{i},y_{-i}^{0},y_{i},u_{i})\nsubseteq
C_{i}(x_{-i}^{0})$ for each $u_{i}\in \overline{A}_{i}(x_{-i}^{0}).$

Let $(x_{i}^{1},y_{i}^{1}),(x_{i}^{2},y_{i}^{2})\in
D(x_{-i}^{0},y_{-i}^{0}). $

This means $f_{i}(x_{-i}^{0},x_{i}^{1},y_{-i}^{0},y_{i}^{1},u_{i})\nsubseteq
C_{i}(x_{-i}^{0})$ for each $u_{i}\in \overline{A}_{i}(x_{-i}^{0})$ and $%
f_{i}(x_{-i}^{0},x_{i}^{2},y_{-i}^{0},y_{i}^{2},u_{i})\nsubseteq
C_{i}(x_{-i}^{0})$ for each $u_{i}\in \overline{A}_{i}(x_{-i}^{0}).$

Let $(x_{i}(\lambda ),y_{i}(\lambda ))=\lambda
(x_{i}^{1},y_{i}^{1})+(1-\lambda )(x_{i}^{2},y_{i}^{2})$ be defined for each 
$\lambda \in \lbrack 0,1].$

We claim that $(x_{i}(\lambda ),y_{i}(\lambda ))\in D(x_{-i}^{0},y_{-i}^{0})$
for each $\lambda \in \lbrack 0,1].$

Suppose, on the contrary, that there exist $\lambda _{0}\in \lbrack 0,1]$ $\ 
$\ and $u_{i}(\lambda _{0})$ such that $f_{i}(x_{-i}^{0},x_{i}(\lambda
_{0}),y_{-i}^{0},y_{i}(\lambda _{0}),u_{i}(\lambda _{0}))\subseteq
C_{i}(x_{-i}^{0}).$ Since $f_{i}(x_{-i}^{0},\cdot ,y_{-i}^{0},\cdot
,u_{i}(\lambda _{0})):X_{i}\times X_{i}\rightarrow 2^{Y_{i}}$ is\textit{\ }$%
C_{i}(x_{-i}^{0})-$\textit{\ }quasi-convex$,$ we have that:

$f_{i}(x_{-i}^{0},x_{i}^{1},y_{-i}^{0},y_{i}^{1},u_{i}(\lambda
_{0}))\subseteq f_{i}(x_{-i}^{0},x_{i}(\lambda ),y_{-i}^{0},y_{i}(\lambda
),u_{i}(\lambda _{0}))+C_{i}(x_{-i}^{0})$ or

$f_{i}(x_{-i}^{0},x_{i}^{2},y_{-i}^{0},y_{i}^{2},u_{i}(\lambda
_{0}))\subseteq $ $f_{i}(x_{-i}^{0},x_{i}(\lambda ),y_{-i}^{0},y_{i}(\lambda
),u_{i}(\lambda _{0}))+C_{i}(x_{-i}^{0}).$\newline
On the other hand, it is true that $f_{i}(x_{-i}^{0},x_{i}(\lambda
_{0}),y_{-i}^{0},y_{i}(\lambda _{0}),u_{i}(\lambda _{0}))\subseteq
C_{i}(x_{-i}^{0}).$ We obtain that:

$f_{i}(x_{-i}^{0},x_{i}^{j},y_{-i}^{0},y_{i}^{j},u_{i}(\lambda
_{0}))\subseteq C_{i}(x_{-i}^{0})+C_{i}(x_{-i}^{0})\subseteq $

$\subseteq C_{i}(x_{-i}^{0})$ for $j=1$ or for $j=2$.

This contradicts the assumption that $%
(x_{i}^{1},y_{i}^{1}),(x_{i}^{2},y_{i}^{2})\in D(x_{-i}^{0},y_{-i}^{0})$.
Consequently, $D(x_{-i}^{0},y_{-i}^{0})$ is convex.

All assumptions of Theorem 4 are verified and there exists $(x^{\ast
},y^{\ast })$ an equilibrium for the associated generalized abstract economy 
$(X_{i},A_{i},F_{i},P_{i})_{i\in N}$. Then, for each $i\in N,$ $x_{i}^{\ast
}\in A_{i}(x_{-i}^{\ast }),$ $y_{i}^{\ast }\in F_{i}(y_{-i}^{\ast })$ and $%
f_{i}(x^{\ast },y^{\ast },u_{i})\nsubseteq $int$C_{i}(x_{-i}^{\ast })$ for
each $u_{i}\in A_{i}(x_{-i}^{\ast }),$ that is, $(x^{\ast },y^{\ast })$ is a
solution for SGVQEP (IV).\medskip 

\textit{Proof of Theorem 13. }For each $i\in N,$ let us define $%
P_{i}:X\times X\rightarrow 2^{X_{i}},$ by

$P_{i}(x,y)=\{u_{i}\in X_{i}:$ $f_{i}(x,y,u_{i})\subseteq $int$%
C_{i}(x_{-i})\}$ for each $(x,y)\in X\times X.$

Assumption d) implies that $\{(x_{i},y_{i})\in X_{i}\times
X_{i}:A_{i}(x_{-i})\cap P_{i}(x,y)=\emptyset \}\mathit{\ }$ is nonempty for
each $(x_{-i},y_{-i})\in X_{-i}\times X_{-i}.$

Let us consider $(x_{-i}^{0},y_{-i}^{0})\in X_{-i}\times X_{-i}$ and let $%
D_{i}(x_{-i}^{0},y_{-i}^{0})=\{(x_{i},y_{i})\in X_{i}\times
X_{i}:A_{i}(x_{-i}^{0})\cap P_{i}(x_{-i}^{0},x_{i},y_{-i}^{0},y_{i}).$ We
can prove that $D_{i}(x_{-i}^{0},y_{-i}^{0})$ is convex, as in the proof of
Theorem 12.

By using a similar argument as in the proof of Theorem 9, we can prove that
Gr$P_{i}$ is open.

All assumptions of Theorem 5 are verified. According to this result, there
exists $(x^{\ast },y^{\ast })$ an equilibrium for the associated generalized
abstract economy $(X_{i},A_{i},F_{i},P_{i})_{i\in N}.$ Then, for each $i\in
N,$ $x_{i}^{\ast }\in \overline{A}_{i}(x_{-i}^{\ast }),$ $y_{i}^{\ast }\in 
\overline{F}_{i}(y_{-i}^{\ast })$ and $f_{i}(x^{\ast },y^{\ast
},u_{i})\nsubseteq $int$C_{i}(x_{-i}^{\ast })=\emptyset $ for each $u_{i}\in
A_{i}(x_{-i}^{\ast }),$ that is, $(x^{\ast },y^{\ast })$ is a solution for
SGVQEP (IV).

\end{document}